\documentclass[11pt]{article}
\usepackage[T1]{fontenc}
\usepackage{lmodern,amsmath,amssymb,amsthm,mathtools,mathrsfs}
\usepackage[a4paper,margin=28mm]{geometry}
\usepackage{microtype}
\usepackage{needspace}
\usepackage[hidelinks]{hyperref}
\numberwithin{equation}{section}
\newtheorem{theorem}{Theorem}[section]
\newtheorem{proposition}[theorem]{Proposition}
\newtheorem{lemma}[theorem]{Lemma}
\newtheorem{corollary}[theorem]{Corollary}
\DeclareMathOperator{\Var}{Var}
\DeclareMathOperator{\tr}{tr}
\DeclareMathOperator{\Hess}{Hess}
\DeclareMathOperator{\Sec}{Sec}

\newcommand{\R}{\mathbb R}
\newcommand{\dd}{\,d}
\newcommand{\eps}{\varepsilon}
\allowdisplaybreaks[2]
\title{Logarithmic Sobolev and exponential inequalities
for submanifolds in Cartan--Hadamard manifolds}
\author{%
Hongyu Chen\thanks{School of Mathematics, Sichuan University,
Chengdu 610064, China.
Email: \texttt{hongyu.chern@gmail.com}.}
\qquad
Li Ma\thanks{Corresponding author.
School of Mathematics and Physics,
University of Science and Technology Beijing,
Beijing 100083, China.
Email: \texttt{nuslma@126.com}.}
}
\date{August 2026}

\begin{document}
\maketitle

\begin{abstract}
We establish a logarithmic Sobolev inequality for complete $C^{1,1}$ immersed submanifolds of dimension $m\ge2$ in Cartan--Hadamard manifolds with sectional curvature at most $-b^2$. The inequality retains the Euclidean sharp constant after subtraction of $(m-1)^2b^2$ from the Fisher and mean-curvature energy; this coefficient is maximal for a lower bound uniform over $b\ge0$. For immersed surfaces, we prove both a finite-domain critical exponential inequality and a global Ruf-type inequality with a positive $L^2$ term. Both have potential $|H|^2-b^2$, exponent $16\pi$ and an absolute multiplicative constant, and follow from one geometric heat-kernel comparison. The exponent and curvature coefficient are optimal. We classify Gaussian equality in the logarithmic inequality and geometric equality in the higher-dimensional Nash consequence. All inequality constants are independent of the codimension.

\medskip
\emph{AMS Classification 2020:} Primary 46E35; Secondary 53C42, 58J35.

\end{abstract}

\tableofcontents{}

\section{Introduction}

Log-Sobolev inequalities on complete Riemannian manifolds are basic analytic tools in geometric analysis. 
In the study of Ricci flow theory, G.Perelman \cite{P02} introduced W-functional to overcome the difficult collapsing problem 
during the evolution of Riemannian metrics. 
The W-functional is a kind of log-Sobolev inequality about the Ricci flow and this has been further developed by L. Ni \cite{N04}  
on a fixed complete Riemannian manifold. The question that whether the infimum of related W-functional can be obtainable has been 
considered in the paper \cite{M20}. In recent preprint \cite{OAI26}, 
OpenAI claims the proof of generalized Cartan-Hadamard isoperimetry and related rigidity. 
A natural equation is whether the log-Sobolev inequality is still true for immersed submanifolds of a Riemannian manifold. 
In a Cartan--Hadamard ambient manifold, one of our aims in this paper is to find the optimal constants and curvature corrections in the log-Sobolev inequality, together with the geometry of equality.

Mean curvature compensates for the geometry of an immersed submanifold in Sobolev inequalities, as in Michael and Simon \cite{MS73} and Hoffman and Spruck \cite{HS74,HS75}. Sharp power inequalities now provide a starting point for studying entropy and exponential integrability. The questions considered here are whether the Euclidean logarithmic constant persists after subtraction of the critical negative-curvature term, and whether the corresponding finite-domain and global surface estimates hold with multiplicative constants uniform over the immersed geometry.

Fix integers $m\ge2$ and $k\ge1$. Throughout the paper, $(M^{m+k},g)$ is smooth, complete and simply connected, with $\Sec_M\le-b^2$ for a parameter $b\ge0$ unless only nonpositive curvature is specified, and $\iota:\Sigma^m\to M$ is a boundaryless $C^{1,1}$ immersion whose induced length metric is complete. The second fundamental form $II$ and the mean curvature vector $H$ are defined almost everywhere. We use the unaveraged convention $H=\operatorname{tr}_{g_\Sigma}II$, where $g_\Sigma$ is the induced metric, as in Ma and Yang \cite{MY04} and Brendle \cite{Br22}. All integrals are source integrals with respect to the induced measure $d\sigma$, including immersion multiplicities. The source may be disconnected. Gradients are weak gradients and $\Delta=\operatorname{div}\nabla$. The codimension $k$ is fixed throughout each variational argument. We write $|\mathbb S^m|$ for the volume of the Euclidean unit $m$-sphere.

\paragraph{The logarithmic curvature correction.}
Ecker \cite{Ec00} proved logarithmic Sobolev inequalities for Euclidean submanifolds. Brendle \cite{Br22} obtained the sharp compact inequality, and Balogh and Krist\'aly \cite{BK26} established the complete Euclidean inequality and its Gaussian rigidity. Yi and Zheng \cite{YZ24} considered compact submanifolds in nonnegatively curved ambient manifolds under an asymptotic Gaussian-volume hypothesis. Our first main result concerns the subtraction $(m-1)^2b^2$ under an ambient sectional-curvature upper bound.

\Needspace{10\baselineskip}
\begin{theorem}\label{thm:main}
Suppose $\Sec_M\le-b^2$ for some $b\ge0$. Let $f\in W^{1,1}_{\mathrm{loc}}(\Sigma)$ be positive almost everywhere, with $\int_\Sigma f\dd\sigma=1$ and
\begin{equation}\label{eq:finite}
\int_\Sigma\frac{|\nabla f|^2}{f}\dd\sigma<\infty,
\qquad
\int_\Sigma\bigl(|H|^2f+f|\log f|\bigr)\dd\sigma<\infty.
\end{equation}
For every $0\le\beta\le(m-1)^2$,
\begin{equation}\label{eq:main}
\int_\Sigma\left(\frac{|\nabla f|^2}{f}-f\log f+
(|H|^2-\beta b^2)f\right)\dd\sigma
\ge m+\frac m2\log(4\pi).
\end{equation}
The constant on the right is independent of $k$ and optimal over all the stated data for each fixed $k$. Equality holds if and only if $b=0$ and there is an isometry $\Phi:\Sigma\to\R^m$ for which $\iota\circ\Phi^{-1}$ is a complete totally geodesic isometric embedding and
$$
f\circ\Phi^{-1}(x)=(4\pi)^{-m/2}e^{-|x-x_0|^2/4},
$$
for some $x_0\in\R^m$.
\end{theorem}

The coefficient $(m-1)^2$ is maximal when the lower bound is required to be uniform over $b\ge0$: every larger coefficient makes the functional unbounded below over the stated class. For each fixed $b$, subtracting a larger coefficient still gives the explicit lower bound in Proposition~\ref{prop:threshold}; the best lower bound at a prescribed $b>0$ is a separate question. Ambient homotheties give its homogeneous form in Corollary~\ref{cor:homogeneous}, whose Gaussian extremizers have arbitrary scale.

For $m\ge3$, Do, Ganguly, Lam and Lu \cite[Theorem~1.3]{DGLL26} proved an intrinsic logarithmic inequality on $\mathbb H^m(-1)$ with the same Euclidean constant and subtraction $m^2(m-1)/(m+2)$ in the Fisher normalization. On a totally geodesic $\mathbb H^m(-1)$, Theorem~\ref{thm:main} and Corollary~\ref{cor:homogeneous} instead subtract $(m-1)^2$. The difference is
$$
(m-1)^2-\frac{m^2(m-1)}{m+2}=\frac{(m-1)(m-2)}{m+2}>0.
$$
Here $H=0$, so this difference is independent of the mean-curvature convention.

\paragraph{Critical exponential inequalities on surfaces.}
Our second main result concerns exponential integrability under a bound on the energy containing $4|\nabla v|^2$ and the potential $|H|^2-b^2$. The critical exponent $16\pi$ corresponds to the usual exponent $4\pi$ after the substitution $u=2v$. The geometric content is the potential and the uniform multiplicative constant. We distinguish the finite-domain estimate, whose right-hand side is proportional to area, from the global estimate, which remains meaningful on sources of infinite area.

\Needspace{29\baselineskip}
\begin{theorem}\label{thm:critical-exponential}
Suppose $m=2$ and $\Sec_M\le-b^2$. There is an absolute finite constant $\mathcal C$ with the following properties.
\begin{enumerate}
\item[(i)] For every relatively compact smooth source domain $U\subset\Sigma$ of positive area and every real $v\in W^{1,2}_0(U)$,
\begin{equation}\label{eq:critical-exponential}
\int_U[4|\nabla v|^2+(|H|^2-b^2)v^2]\dd\sigma\le1
\quad\Longrightarrow\quad
\int_U e^{16\pi v^2}\dd\sigma\le\mathcal C|U|.
\end{equation}
When $\Sigma$ is compact, one may take $U=\Sigma$ and $v\in W^{1,2}(\Sigma)$.
\item[(ii)] For every $\lambda>0$ and every real $v\in W^{1,2}(\Sigma)$ with $|H|v\in L^2(\Sigma)$,
\begin{equation}\label{eq:global-constraint}
\int_\Sigma[4|\nabla v|^2+(|H|^2-b^2+\lambda)v^2]\dd\sigma\le1
\end{equation}
implies
\begin{equation}\label{eq:global-exponential}
\lambda\int_\Sigma(e^{16\pi v^2}-1)\dd\sigma\le\mathcal C.
\end{equation}
\end{enumerate}
The constant is independent of the geometry, codimension, $b$, $U$ and $\lambda$. The exponent $16\pi$ is optimal in both parts, in part~(ii) even for each fixed $\lambda>0$. The curvature coefficient one is maximal in both uniform statements: replacing $-b^2$ by $-\beta b^2$ with $\beta>1$ invalidates part~(i) already on one fixed totally geodesic hyperbolic plane, and invalidates part~(ii), for each prescribed $\lambda>0$, on a suitably chosen totally geodesic hyperbolic plane. The best multiplicative constants and equality cases are not determined here.
\end{theorem}

Part~(ii) is a full-energy, or Ruf-type, endpoint. On a Euclidean plane it reduces, after $u=2v$ and scaling, to the classical full-norm inequality recalled and extended by Li and Ruf \cite{LR08}. The subtraction of one in the exponential is necessary on infinite-area sources. The positive parameter $\lambda$ supplies the time scale and the $L^2$ control; finite-domain exhaustion alone does not yield \eqref{eq:global-exponential}. The theorem allows self-intersections and nonproper immersions, with all masses computed on the source.

Hyperbolic Hardy--Moser inequalities already involve subtraction of the bottom of the spectrum. Lu and Yang \cite[Theorem~1.2 and Section~4]{LY16} give a proof of the strengthened hyperbolic inequality attributed there to Mancini, Sandeep and Tintarev, with the regularized exponential integrated over the whole hyperbolic plane. Morpurgo and Qin \cite[Theorem~1]{MQ24} prove intrinsic inequalities on Hadamard surfaces with two-sided negative-curvature bounds, using the Dirichlet energy. Do, Lam, Lu and Ponciano \cite{DLLP26} obtain submanifold exponential inequalities under nonnegative ambient curvature and a small total mean-curvature assumption. Theorem~\ref{thm:critical-exponential} concerns the displayed mean-curvature potential, an ambient curvature upper bound, and both estimates uniform over the geometry; its global part also allows infinite area and a signed potential without an assumed spectral gap.

Bernstein and Bhattacharya \cite[Section~5]{BB25} develop Gaussian and hyperbolic entropy functionals for mean curvature flow in Cartan--Hadamard manifolds. Their local kernel calculation is closely related to the normal-square identity used here. We use an explicit radial supersolution, so no model heat-kernel convexity theorem is needed. Source-sheet initial mass, the global Friedrichs form and Dirichlet exhaustion give one diagonal estimate for $-4\Delta+|H|^2-b^2$. Separate spectral arguments then give the two critical endpoints.

\paragraph{Power Sobolev estimates and their consequences.}
The higher-dimensional power inequality is an input to our logarithmic argument. Case and Mengesha \cite[Theorems~1.1 and~1.3, Corollary~1.4]{CM26} prove the sharp $L^2$ Michael--Simon inequality in Euclidean space and simply connected spaceforms, in a conformally invariant formulation for conformally developable ambient manifolds. Their formula~(1.7) also gives the curvature-bounded estimate in that setting after changing from the averaged mean curvature $H/m$ to our convention. Balogh, Krist\'aly and Mester \cite{BKM26} obtain explicit codimension-independent constants for Euclidean minimal submanifolds, asymptotically sharp in the dimension.

For general Cartan--Hadamard ambient manifolds, the preprint \cite[Theorem~3.1]{OAI26} already states the exact power inequality for locally $C^{1,1}$ immersed hypersurfaces, allowing incomplete sources and boundary when the tests are supported in the interior. Its normalization uses $|H|^2=m^2h^2$; multiplication by $m$, followed by scaling from $b=1$ when $b>0$, gives \eqref{eq:sharp-sobolev} in codimension one. Its radial comparison maps and centered coordinate variations extend to fixed arbitrary codimension by replacing the scalar normal square by the full normal-vector square. Section~\ref{sec:power} records the normal-vector calculation needed for this extension and the comparison identities used in the later equality analysis.

Applying the power inequality to product immersions and passing to the product limit gives the sharp unshifted logarithmic inequality in every dimension $m\ge2$; Proposition~\ref{prop:unshifted} supplies this short derivation. When $b>0$, product ambient manifolds have mixed zero-curvature planes, so this argument retains nonpositive curvature but loses the strict upper bound $-b^2$. The proof of the critical subtraction in Theorem~\ref{thm:main} therefore requires a separate argument. The higher-dimensional Nash inequality is a H\"older consequence of the power estimate, and the surface Nash estimate follows from the logarithmic inequality. They are stated as Corollary~\ref{thm:nash}; the higher-dimensional equality proof recovers an entire totally geodesic model ball in the ambient manifold. Its optimal constant over varying submanifolds is distinct from the sharp Nash constant on a fixed Euclidean space studied by Carlen and Loss \cite{CL93}.

\paragraph{Proofs and organization.}
For the logarithmic correction, a sequence approaching the infimum satisfies an approximate weighted Poincar\'e inequality by normalized second variation. A bounded radial vector field permits centering without a prescribed moment condition. Hessian and Rauch comparison yield a positive remainder on each fixed annulus. In dimensions $m\ge3$, a Hardy estimate prevents concentration at the center; in dimension two a local divergence estimate serves this role. Escape to infinity then contradicts the Poincar\'e estimate unless the curvature parameters tend to zero. The unshifted inequality gives the required lower bound. The Poincar\'e estimate used in this argument is a property of the chosen approximate minimizers, not an assumption on arbitrary densities.

The plan of the paper is below. Section~\ref{sec:prelim} records the weak geometric identities and the density convention. Section~\ref{sec:power} recalls the power inequality, explains its codimension extension, and derives the unshifted logarithmic inequality. Section~\ref{sec:analytic} supplies coercivity, Hardy control and moment estimates. Sections~\ref{sec:shift} and~\ref{sec:equality} prove the logarithmic correction and its equality classification for $m\ge3$, and Section~\ref{sec:log-two} treats dimension two. Section~\ref{sec:sharpness} proves the logarithmic threshold and homogeneous form. Section~\ref{sec:critical-two} proves the common heat-kernel comparison and both exponential inequalities, and Section~\ref{sec:nash} gives the Nash consequences and their geometric sharpness.

\section{Weak geometric identities}\label{sec:prelim}

The induced metric has locally Lipschitz coefficients, and $H$ is locally essentially bounded. For a smooth ambient function $\psi$ the weak chain rule is
\begin{equation}\label{eq:chain}
\Delta_\Sigma(\psi\circ\iota)=\tr_{T\Sigma}\Hess^M\psi+\langle H,\nabla^M\psi\rangle.
\end{equation}
For a nonnegative density $f$ with $\sqrt f\in W^{1,2}(\Sigma)$ and finite $\int|H|^2f$, a smooth ambient vector field $Y$ satisfies the following identity when its displayed terms are integrable and its cutoff boundary term tends to zero:
\begin{equation}\label{eq:weighted-div}
\int\tr_{T\Sigma}\nabla^MY\,f\dd\sigma
=-\int\left\langle\frac{\nabla f}{f}+H,Y\right\rangle f\dd\sigma.
\end{equation}
The quotient is used only on $\{f>0\}$. The gradient of a Sobolev function vanishes almost everywhere on its zero set, so the vector in \eqref{eq:weighted-div} has squared $L^2(fd\sigma)$ norm $\int_\Sigma(\frac{|\nabla f|^2}{f}+|H|^2f)\dd\sigma$.

These identities are first paired with compact source tests. On a complete source, intrinsic Lipschitz cutoffs have compact support and gradients bounded by $C/R$. If the source is disconnected, exhaust it by finitely many components and use intrinsic cutoffs on each of them; integrability makes the discarded component tails tend to zero. For a bounded ambient field their boundary error tends to zero. For a field bounded by $C(1+r)$, where $r$ is ambient distance from a fixed point, it tends to zero if $\int r^2f<\infty$. If a Hessian trace is nonnegative, Fatou's lemma first proves its integrability from the finite right side; the full identity then follows. This order will be used below.

In the logarithmic arguments the variational unknowns are probability densities $f$. Finite energy means $\sqrt f\in W^{1,2}(\Sigma)$ and $\int_\Sigma|H|^2f\dd\sigma<\infty$; finite absolute entropy means $\int_\Sigma f|\log f|\dd\sigma<\infty$. For nonnegative densities, we set $f\log f=0$ and $|\nabla f|^2/f=0$ on $\{f=0\}$. Constants with indicated subscripts depend only on those subscripts and may change between estimates.

\begin{lemma}\label{lem:density}
Let $f\in W^{1,1}_{\mathrm{loc}}(\Sigma)$ be nonnegative with $\int_\Sigma f\dd\sigma=1$. Then the conditions in \eqref{eq:finite} are equivalent to $\sqrt f\in W^{1,2}(\Sigma)$, $\int_\Sigma|H|^2f\dd\sigma<\infty$ and $\int_\Sigma f|\log f|\dd\sigma<\infty$. In particular,
$$
\int_\Sigma\frac{|\nabla f|^2}{f}\dd\sigma=4\int_\Sigma|\nabla\sqrt f|^2\dd\sigma.
$$
\end{lemma}
\begin{proof}
Apply the local Sobolev chain rule to $\sqrt{f+\eps}$. Its gradient squared is $|\nabla f|^2/[4(f+\eps)]$, bounded in integral by one quarter of the Fisher information. Local weak compactness and exhaustion give $\sqrt f\in W^{1,2}(\Sigma)$, whose squared $L^2$ norm is one. The chain rule on $\{f>\eps\}$, followed by $\eps\downarrow0$, identifies its gradient on $\{f>0\}$. The gradient of $\sqrt f$ vanishes almost everywhere on $\{f=0\}$, which proves the identity. Conversely, the product rule for $f=(\sqrt f)^2$ gives $f\in W^{1,1}(\Sigma)$ and $\nabla f=2\sqrt f\,\nabla\sqrt f$. The curvature and entropy conditions are unchanged.
\end{proof}

\section{The power inequality and the unshifted logarithmic bound}\label{sec:power}

We use \cite[Theorem~3.1]{OAI26} in our mean-curvature normalization. For arbitrary fixed codimension, only the normal square and tangent trace in that proof need to be changed. We give those interfaces below and refer to the source for the common compactness and centering arguments.

\begin{proposition}\label{thm:sharp-sobolev}
Suppose $m\ge3$ and $\Sec_M\le-b^2$. Then for every real $v\in C_c^1(\Sigma)$,
\begin{equation}\label{eq:sharp-sobolev}
\int_\Sigma\left[\frac{4m}{m-2}|\nabla v|^2+(|H|^2-m^2b^2)v^2\right]\dd\sigma
\ge m^2|\mathbb S^m|^{2/m}\|v\|_{2m/(m-2)}^2.
\end{equation}
For every $0\le\beta\le m^2$, the same constant holds with $m^2b^2$ replaced by $\beta b^2$. The same statement holds for compactly supported $W^{1,2}$ tests. Its constant is optimal uniformly over the stated class for each fixed $k$ and is independent of $k$.
\end{proposition}
\begin{proof}[The codimension extension]
The codimension-one assertion is \cite[Theorem~3.1]{OAI26}, after multiplication by $m$ and scaling. For general $k$, let $F:M\to\mathbb S^{m+k}$ and $\Omega>0$ be the radial comparison map and weight of \cite[Lemma~2.1]{OAI26}, expressed for general $b$ by scaling. They satisfy
\begin{equation}\label{eq:comparison-interface}
F^*g_{\mathbb S^{m+k}}\le\Omega^2g,\qquad
\Hess^M\log\Omega\le d\log\Omega\otimes d\log\Omega
-\tfrac12(b^2+\Omega^2+|\nabla^M\log\Omega|^2)g.
\end{equation}
The weak restriction identity \eqref{eq:chain}, followed by integration by parts and completion of the full normal-vector square, gives
\begin{align}\label{eq:vector-square}
&\int_\Sigma\left[\frac{4m}{m-2}|\nabla v|^2+(|H|^2-m^2b^2)v^2\right]\dd\sigma\notag\\
&\quad\ge m^2\int_\Sigma\Omega^2v^2\dd\sigma
+\int_\Sigma|H-m(\nabla^M\log\Omega)^\perp|^2v^2\dd\sigma\notag\\
&\qquad+\frac{4m}{m-2}\int_\Sigma
\left|\nabla v-\frac{m-2}{2}v\nabla_\Sigma\log\Omega\right|^2\dd\sigma.
\end{align}
Indeed, the second inequality in \eqref{eq:comparison-interface} implies
$$
\begin{aligned}
|H|^2-m^2b^2\ge{}&m^2\Omega^2+|H-m(\nabla^M\log\Omega)^\perp|^2\\
&+m(m-2)|\nabla_\Sigma\log\Omega|^2+2m\Delta_\Sigma\log\Omega.
\end{aligned}
$$
Multiplication by $v^2$ yields \eqref{eq:vector-square}. The coordinate sum satisfies
\begin{equation}\label{eq:short}
\sum_{\alpha=1}^{m+k+1}|\nabla_\Sigma(F_\alpha\circ\iota)|^2\le m\Omega^2.
\end{equation}
Thus the normal term contains every normal direction in one square, and the coordinate trace depends on $m$, not on the number of coordinates.

The local compactness lemmas \cite[Lemmas~3.3--3.4]{OAI26} require a locally Lipschitz induced metric and a locally bounded scalar potential. They apply here to $|H|^2-m^2b^2$. The pushforward of a compactly supported critical mass is nonatomic, since an immersion fiber is discrete; \cite[Lemma~2.3]{OAI26} therefore centers the $m+k+1$ coordinates. The second-variation calculation in \cite[Section~3.3]{OAI26}, with \eqref{eq:short}, gives the same threshold. Formula~\eqref{eq:vector-square} forces the quotient minimizer divided by $\Omega^{(m-2)/2}$ to be constant. The zero-extension contradiction and the point cutoff on a compact component are exactly those in that section; their only dimensional restriction is $m>2$. This proves the stated extension without any codimension factor. Adding $(m^2-\beta)b^2\|v\|_2^2$ gives the remaining coefficients. Round Euclidean spheres give the uniform sharpness.
\end{proof}

Round Euclidean spheres with constant tests attain the displayed constant. The coefficient $m^2$ is maximal even at each fixed $b>0$: for every larger coefficient, constant tests on sufficiently large hyperbolic geodesic spheres have negative energy, as computed in Corollary~\ref{cor:nash-range}. These sharpness facts belong to the power estimate and its model tests.

\begin{proposition}\label{prop:unshifted}
Let $m\ge2$ and $\Sec_M\le0$. Every nonnegative probability density $f$ of finite energy and finite absolute entropy satisfies
\begin{equation}\label{eq:unshifted-homogeneous}
\int_\Sigma\left(\frac{|\nabla f|^2}{f}+|H|^2f\right)\dd\sigma
\ge2\pi e m\exp\left(\frac2m\int_\Sigma f\log f\dd\sigma\right).
\end{equation}
Consequently,
$$
\int_\Sigma\left(\frac{|\nabla f|^2}{f}+|H|^2f-f\log f\right)\dd\sigma
\ge m+\frac m2\log(4\pi).
$$
\end{proposition}
\begin{proof}
We use the product-limit argument of Beckner and Pearson \cite{BP98}, applied to the power inequality for product immersions.
First suppose that $f$ is bounded and compactly supported and $\sqrt f\in W^{1,2}(\Sigma)$. For an integer $q$ with $qm\ge3$, use the product immersion $\Sigma^q\to M^q$ and the test
$$
V_q(x_1,\ldots,x_q)=\prod_{l=1}^q\sqrt{f(x_l)}.
$$
The product ambient manifold is complete, simply connected and nonpositively curved. The immersion is complete and $C^{1,1}$, and has finite codimension $qk$. Its mean-curvature vector has squared norm $\sum_{l=1}^q|H(x_l)|^2$. Since $\int f=1$, Proposition~\ref{thm:sharp-sobolev} gives
\begin{align*}
q\left[\frac{qm}{qm-2}\int_\Sigma\frac{|\nabla f|^2}{f}\dd\sigma
+\int_\Sigma|H|^2f\dd\sigma\right]
\ge (qm)^2|\mathbb S^{qm}|^{2/(qm)}
\left(\int_\Sigma f^{qm/(qm-2)}\dd\sigma\right)^{(qm-2)/m}.
\end{align*}
The constant is independent of the product codimension. Boundedness and compact support justify differentiation of $\int f^{1+h}$ at $h=0$, giving
$$
\frac{qm-2}{m}\log\int_\Sigma f^{qm/(qm-2)}\dd\sigma
\longrightarrow\frac2m\int_\Sigma f\log f\dd\sigma.
$$
Stirling's formula and $|\mathbb S^n|=2\pi^{(n+1)/2}/\Gamma((n+1)/2)$ give
$$
\lim_{n\to\infty}n|\mathbb S^n|^{2/n}=2\pi e.
$$
Divide the preceding power inequality by $q$ and let $q\to\infty$ to obtain \eqref{eq:unshifted-homogeneous} in this class.

For a general $f$, choose complete intrinsic cutoffs, including finitely many components at each stage, and truncate the values of $\sqrt f$. Normalize the squared tests to have integral one. They are bounded and compactly supported. Their Fisher and mean-curvature energies converge to those of $f$ by $W^{1,2}$ convergence and dominated convergence. Their entropies also converge: the value tails are controlled by $f|\log f|$, and the extra cutoff term is dominated by $f/e$, since $|s\log s|\le1/e$ for $0\le s\le1$. Normalization factors tend to one. Passing to the limit proves \eqref{eq:unshifted-homogeneous} for every stated density.

Its energy is positive. To deduce the additive form, minimize
$$
x-\frac m2\log\frac{x}{2\pi e m},\qquad x>0.
$$
The minimum occurs at $x=m/2$ and equals $m+(m/2)\log(4\pi)$.
\end{proof}

The product argument gives the unshifted bound also for $m=2$. It does not supply the negative-curvature correction or an equality classification.

\section{Coercivity, Hardy estimates and moments}\label{sec:analytic}

\begin{lemma}\label{lem:sobolev}
For $m\ge3$, on every stated immersion with $\Sec_M\le0$,
$$
\|v\|_{2m/(m-2)}^2
\le\frac{\int_\Sigma(4|\nabla v|^2+|H|^2v^2)\dd\sigma}
{m(m-2)|\mathbb S^m|^{2/m}},
$$
for every $v\in W^{1,2}(\Sigma)$ with $|H|v\in L^2$. For $\|v\|_2=1$ and finite absolute entropy,
\begin{align*}
\int_\Sigma v^2\log v^2\dd\sigma
&\le\frac m2\log\left\{
\frac{\int_\Sigma(4|\nabla v|^2+|H|^2v^2)\dd\sigma}
{m(m-2)|\mathbb S^m|^{2/m}}\right\},\\
\int_\Sigma(4|\nabla v|^2+|H|^2v^2-v^2\log v^2)\dd\sigma
&\ge\frac m2\left[1+\log\bigl(2(m-2)|\mathbb S^m|^{2/m}\bigr)\right].
\end{align*}
\end{lemma}
\begin{proof}
Proposition~\ref{thm:sharp-sobolev} with $b=0$ gives
$$
m^2|\mathbb S^m|^{2/m}\|v\|_{2m/(m-2)}^2
\le\frac{4m}{m-2}\int|\nabla v|^2+\int|H|^2v^2
\le\frac m{m-2}\int(4|\nabla v|^2+|H|^2v^2),
$$
for compact tests. Complete intrinsic cutoffs followed by value truncation extend it to the displayed finite-energy class: the cutoff gradient error is bounded by $C R^{-2}\|v\|_2^2$, and the curvature terms converge by dominated convergence. Jensen's inequality for the probability $v^2d\sigma$ proves the entropy bound. Minimizing
$$
x-\frac m2\log\frac{x}{m(m-2)|\mathbb S^m|^{2/m}},\qquad x>0,
$$
gives the second bound.
\end{proof}

\begin{lemma}\label{lem:hardy}
If $m\ge3$, $\Sec_M\le-b^2$, $b\ge0$, and $r=d_M(p,\iota(\cdot))$, then for every real $v\in W^{1,2}(\Sigma)$ with $|H|v\in L^2$,
$$
\int\bigl(4|\nabla v|^2+|H|^2v^2\bigr)\dd\sigma-(m-1)^2b^2\|v\|_2^2\ge\int_\Sigma\frac{v^2}{r^2}\dd\sigma.
$$
\end{lemma}
\begin{proof}
Use continuous limits at $b=0$. Hessian comparison and \eqref{eq:chain} give
$$
\Delta_\Sigma r\ge b\coth(br)\bigl(m-1+|(\nabla^Mr)^\perp|^2\bigr)
+\langle H,\nabla^Mr\rangle.
$$
For compact tests away from $r=0$, multiply by
$2((m-1)b\coth(br)-r^{-1})v^2$ and integrate by parts. Completing the tangential and normal squares gives
\begin{align*}
&\int\bigl[4|\nabla v|^2+(|H|^2-(m-1)^2b^2)v^2\bigr]\dd\sigma\\
&\quad\ge4\int\left|\nabla v+
\frac{(m-1)b\coth(br)-r^{-1}}2v\nabla r\right|^2\dd\sigma\\
&\qquad+\int\left|H+
\bigl((m-1)b\coth(br)-r^{-1}\bigr)(\nabla^Mr)^\perp\right|^2v^2\dd\sigma\\
&\qquad+2\int\left[
\bigl((m-1)b\coth(br)-r^{-1}\bigr)b\coth(br)
+\frac{(m-1)b^2}{\sinh^2(br)}-\frac1{r^2}
\right]|(\nabla^Mr)^\perp|^2v^2\dd\sigma\\
&\qquad+\int\left[\frac1{r^2}
+\frac{(m-1)(m-3)b^2}{\sinh^2(br)}\right]v^2\dd\sigma.
\end{align*}
The coefficient in the third integral is obtained by differentiating
$b\coth(br)$, whose derivative is $-b^2/\sinh^2(br)$.
Since $b\coth(br)\ge r^{-1}$, the expression in its square brackets is at least $(m-3)/r^2$. All the remainders are therefore nonnegative. At $b=0$, the last weight equals $(m-2)^2/r^2\ge r^{-2}$.

The fiber $\iota^{-1}(p)$ is closed and discrete and meets a compact source set in finitely many points. Removing small coordinate balls around these points has $W^{1,2}$ cutoff cost $O(\eps^{m-2})$ for bounded compact tests. Since $m\ge3$, the cost vanishes. Weak chain rules, local approximation, complete cutoffs and value truncation extend the inequality to every finite-energy $v$.
\end{proof}

\begin{lemma}\label{lem:moments}
Let $f$ be a probability density of finite energy. If
\begin{equation}\label{eq:pi}
\Var_{fd\sigma}(\phi)\le\kappa\int|\nabla\phi|^2f\dd\sigma,
\end{equation}
for every bounded $\phi\in W_{\mathrm{loc}}^{1,2}(\Sigma)$ satisfying $\int_\Sigma|\nabla\phi|^2f\dd\sigma<\infty$, then $d_M(p,\iota(\cdot))\in L^2(fd\sigma)$ for every $p\in M$.
If $\Sec_M\le0$ and $\kappa=2$, there is a center $p$ such that
\begin{align}\label{eq:center-bound}
\int\log_p(\iota(x))f(x)\dd\sigma(x)=0,
\qquad &\int r^2f\dd\sigma\le2m,\notag\\
m\le\int\tr_{T\Sigma}\Hess^M(r^2/2)f\dd\sigma
&\le\sqrt{\int\bigl(\frac{|\nabla f|^2}{f}+|H|^2f\bigr)\dd\sigma\int r^2f\dd\sigma}.
\end{align}
\end{lemma}
\begin{proof}
Apply \eqref{eq:pi} to the bounded truncations $\min(d_M(p_0,\iota(\cdot)),N)$, whose variances are at most $\kappa$. Choose $R_0$ such that
$$
\int_{\{d_M(p_0,\iota(x))\le R_0\}}f\dd\sigma\ge\frac12.
$$
It follows that
$$
\left(\int\min(d_M(p_0,\iota(\cdot)),N)f\dd\sigma-R_0\right)_+^2\le2\kappa,
$$
so the first and second moments of these truncations are bounded independently of $N$. Fatou and the triangle inequality prove the moment assertion.

The function $p\mapsto\int d_M(p,\iota(x))^2f(x)\dd\sigma(x)$ is finite, continuous and proper, and attains its minimum. Its derivative is justified by the first moment, giving the first assertion in \eqref{eq:center-bound}. Rauch comparison gives $|d\log_p|\le1$. Apply \eqref{eq:pi} with $\kappa=2$ to value truncations of its $m+k$ coordinates and pass to the limit in $L^2(fd\sigma)$. If $e_1,\ldots,e_m$ is an orthonormal tangent frame, then
$$
\sum_{\alpha=1}^{m+k}|\nabla(\log_p\circ\iota)_\alpha|^2
=\sum_{i=1}^m|d\log_p(d\iota(e_i))|^2\le m.
$$
Summing the coordinate inequalities therefore gives the second assertion.

Finally $\Hess^M(r^2/2)\ge g$, \eqref{eq:weighted-div}, and Cauchy--Schwarz give the last assertion. The full second moment just established justifies the unbounded radial field.
\end{proof}

\section{The critical curvature correction}\label{sec:shift}

Throughout this section and Section~\ref{sec:equality}, $m\ge3$. The next radial field has divergence $m$ on the model space $\mathbb H^m(-b^2)$. Its boundedness permits centering a probability density without a moment assumption. For $b>0$, its radial coefficient is
$$
\gamma_b(r)=\frac{m}{\sinh^{m-1}(br)}\int_0^r\sinh^{m-1}(bt)\dd t.
$$
It extends smoothly as an odd function across zero and satisfies
\begin{equation}\label{eq:profile-ode}
\gamma_b'+(m-1)b\coth(br)\gamma_b=m.
\end{equation}

\begin{lemma}\label{lem:profile}
For $r>0$,
\begin{align}\label{eq:profile-bounds}
&0<\gamma_b'<1,\qquad 0<\frac{b\gamma_b}{\sinh(br)}<1,\notag\\
&0<b\coth(br)\gamma_b-\gamma_b'<b\gamma_b.
\end{align}
Moreover, $\gamma_b(r)\to m/((m-1)b)$, $\gamma_b'(r)\to0$ and $b\gamma_b(r)/\sinh(br)\to0$ as $r\to\infty$, uniformly for $b$ in compact subintervals of $(0,\infty)$.
\end{lemma}
\begin{proof}
The lower and upper comparison functions $1/(b\coth(br))$ and $m/(mb\coth(br)-b)$ have residuals in \eqref{eq:profile-ode}
$$
-\frac{b^2}{(b\coth(br))^2},\qquad
\frac{m(m+1)b(b\coth(br)-b)}{(mb\coth(br)-b)^2},
$$
respectively. Near zero,
\begin{align*}
\gamma_b(r)&=r-\frac{(m-1)b^2}{3(m+2)}r^3+O(r^5),\\
\frac1{b\coth(br)}&=r-\frac{b^2}3r^3+O(r^5),\qquad
\frac m{mb\coth(br)-b}=r+\frac bmr^2+O(r^3).
\end{align*}
These expansions give $1/(b\coth(br))<\gamma_b<m/(mb\coth(br)-b)$ initially. The first-crossing argument for this linear differential equation therefore gives these strict inequalities on $(0,\infty)$. The upper comparison is smaller than $m/((m-1)b\coth(br))$; the lower comparison and the equation then give $0<\gamma_b'<1$ and $1<b\coth(br)\gamma_b<m/(m-1)$. In particular $\gamma_b<r$, which proves $b\gamma_b/\sinh(br)<1$. Finally
$$
b\coth(br)\gamma_b-\gamma_b'=mb\coth(br)\gamma_b-m>0,
\qquad (mb\coth(br)-b)\gamma_b<m
$$
give the last assertion of \eqref{eq:profile-bounds}. The limit of $\gamma_b$ follows directly from its integral formula by l'H\^opital's rule. The other limits follow from the equation and the exponential denominator. Rescaling $br$ in that formula proves uniformity on every stated parameter interval.
\end{proof}

We record two consequences of these bounds. For $r>0$,
\begin{align}\label{eq:profile-positive}
&b^2\gamma_b^2-\bigl(b\coth(br)\gamma_b-\gamma_b'\bigr)^2>0,\notag\\
&2\bigl(b\coth(br)\gamma_b-\gamma_b'\bigr)
 +(\gamma_b')^2-\frac{b^2\gamma_b^2}{\sinh^2(br)}\ge0.
\end{align}
For the second inequality its left side equals
$$
(\gamma_b'-1)^2+2\bigl(b\coth(br)\gamma_b-1\bigr)
 +1-\frac{b^2\gamma_b^2}{\sinh^2(br)}.
$$
The first expression in \eqref{eq:profile-positive} has a positive lower bound on each compact annulus, uniformly for $b$ in a compact subinterval of $(0,\infty)$.

\begin{lemma}\label{lem:radial}
Let $b>0$, and let $f$ be a probability density of finite energy in the stated geometry. There is a center $p\in M$ such that
\begin{equation}\label{eq:profile-center}
\int\gamma_b(r)\frac{\log_p(\iota(x))}{r}f(x)\dd\sigma(x)=0,
\qquad r=d_M(p,\iota(x)).
\end{equation}
At this center,
\begin{equation}\label{eq:profile-div}
\int\bigl(\frac{|\nabla f|^2}{f}+|H|^2f\bigr)\dd\sigma\int\gamma_b^2f\dd\sigma\ge\left(m+\int(1-|\nabla_\Sigma r|^2)\bigl(b\coth(br)\gamma_b-\gamma_b'\bigr)f\dd\sigma\right)^2.
\end{equation}
If $fd\sigma$ satisfies \eqref{eq:pi} with $\kappa=2$, then
$$
\int\bigl(\frac{|\nabla f|^2}{f}+|H|^2f\bigr)\dd\sigma-(m-1)^2b^2>m/2.
$$
\end{lemma}
\begin{proof}
Fix $p_0\in M$ and minimize
$$
p\longmapsto\int_\Sigma\int_{d_M(p_0,\iota(x))}^{d_M(p,\iota(x))}
\gamma_b(t)\dd t\,f(x)\dd\sigma(x).
$$
The inner integral is bounded in absolute value by $m\,d_M(p,p_0)/((m-1)b)$, so this is finite and continuous and needs no moment hypothesis. Choose $R_0>0$ such that $\int_{\{d_M(p_0,\iota(x))\le R_0\}}f\dd\sigma>1/2$. For sufficiently small $\eta>0$, the limit of $\gamma_b$ gives
$$
\int_0^r\gamma_b(t)\dd t\ge\frac{m(1-\eta)}{(m-1)b}r-C_{b,\eta}.
$$
On the chosen ball this and the triangle inequality give a lower bound with slope $m(1-\eta)/((m-1)b)$ in $d_M(p,p_0)$; on the complement the Lipschitz bound gives slope at least $-m/((m-1)b)$. Choose $\eta>0$ small enough that the resulting slope
$$
\frac{m}{(m-1)b}\left[(2-\eta)\int_{\{d_M(p_0,\iota(x))\le R_0\}}f\dd\sigma-1\right]
$$
is positive. Hopf--Rinow gives a minimizer. Differentiation under the integral is justified by the bounded derivative and yields \eqref{eq:profile-center}.

Rauch comparison bounds the radial and transverse singular values of the differential of $\gamma_b(r)\log_p/r$ by $\gamma_b'$ and $b\gamma_b/\sinh(br)$. Its $m+k$ coordinate functions are bounded, and their summed source gradient square is at most
$$
(\gamma_b')^2|\nabla_\Sigma r|^2+
\frac{b^2\gamma_b^2}{\sinh^2(br)}(m-|\nabla_\Sigma r|^2).
$$
This is the trace on the $m$-dimensional tangent space of the squared differential. Consequently \eqref{eq:pi}, when available, gives
\begin{equation}\label{eq:profile-pi}
\int\gamma_b^2f\dd\sigma\le2\int\left[(\gamma_b')^2|\nabla_\Sigma r|^2+
\frac{b^2\gamma_b^2}{\sinh^2(br)}(m-|\nabla_\Sigma r|^2)\right]f\dd\sigma.
\end{equation}
Hessian comparison and \eqref{eq:profile-ode} give
$$
\tr_{T\Sigma}\nabla^M(\gamma_b\nabla^Mr)
\ge m+(1-|\nabla_\Sigma r|^2)\bigl(b\coth(br)\gamma_b-\gamma_b'\bigr).
$$
Apply \eqref{eq:weighted-div} to this bounded radial field and use Cauchy--Schwarz. Its nonnegative trace is first made integrable by cutoff and Fatou, as in Section~\ref{sec:prelim}. This proves \eqref{eq:profile-div}.

The algebra that combines these two estimates is
\begin{align}\label{eq:radial-algebra}
&\left(m+\int(1-|\nabla_\Sigma r|^2)\bigl(b\coth(br)\gamma_b-\gamma_b'\bigr)f\dd\sigma\right)^2\notag\\
&\quad-(m-1)^2b^2\int\gamma_b^2f\dd\sigma
 -m\int\left[(\gamma_b')^2|\nabla_\Sigma r|^2+
 \frac{b^2\gamma_b^2}{\sinh^2(br)}(m-|\nabla_\Sigma r|^2)\right]f\dd\sigma\notag\\
&\quad=(m-1)\int\left[b^2\gamma_b^2-
 \bigl(b\coth(br)\gamma_b-\gamma_b'\bigr)^2\right]f\dd\sigma\notag\\
&\qquad\quad+\left(\int(1-|\nabla_\Sigma r|^2)\bigl(b\coth(br)\gamma_b-\gamma_b'\bigr)f\dd\sigma\right)^2\notag\\
&\qquad\quad+m\int(1-|\nabla_\Sigma r|^2)\left[2\bigl(b\coth(br)\gamma_b-\gamma_b'\bigr)
 +(\gamma_b')^2-\frac{b^2\gamma_b^2}{\sinh^2(br)}\right]f\dd\sigma.
\end{align}
To verify it, expand the first square in \eqref{eq:radial-algebra}, substitute $m=\gamma_b'+(m-1)b\coth(br)\gamma_b$ in its constant term, and use $b^2\coth^2(br)-b^2=b^2/\sinh^2(br)$. The terms independent of $1-|\nabla_\Sigma r|^2$ reduce to the first integral; the remaining linear terms give the last integral and squaring the integrated remainder gives the remaining quadratic term.

The right side is strictly positive: the first integrand is positive for $r>0$, and the fiber $r=0$ has source measure zero. Combining \eqref{eq:radial-algebra}, \eqref{eq:profile-div} and \eqref{eq:profile-pi} gives
$$
0<\left(\int\bigl(\frac{|\nabla f|^2}{f}+|H|^2f\bigr)\dd\sigma-(m-1)^2b^2-\frac m2\right)\int\gamma_b^2f\dd\sigma.
$$
This proves the strict assertion.
\end{proof}

\begin{proof}[Proof of Theorem~\ref{thm:main}: the inequality for $m\ge3$]
It suffices to prove the critical coefficient $\beta=(m-1)^2$. Fix $B<\infty$ and denote by $a_B$ the infimum of
$$
\int\bigl(\frac{|\nabla f|^2}{f}+|H|^2f\bigr)\dd\sigma-(m-1)^2b^2-\int f\log f\dd\sigma
$$
over the stated positive densities and geometries of the fixed codimension $k$ with $0\le b\le B$. Lemma~\ref{lem:sobolev} bounds $a_B$ from below by a finite constant depending on $m,B$, and the Euclidean Gaussian gives $a_B\le m+(m/2)\log(4\pi)$. Choose data with values at most $a_B+\eps_j$, where $\eps_j\downarrow0$. Their energies and absolute entropies are uniformly bounded, with constants allowed to depend on $B$. For energy this follows from
\begin{align*}
&\int\bigl(\frac{|\nabla f_j|^2}{f_j}+|H_j|^2f_j\bigr)\dd\sigma_j\\
&\quad-\frac m2\log\left\{
\frac{\int(\frac{|\nabla f_j|^2}{f_j}+|H_j|^2f_j)\dd\sigma_j}
{m(m-2)|\mathbb S^m|^{2/m}}\right\}\\
&\quad\le m+\frac m2\log(4\pi)+1+(m-1)^2B^2.
\end{align*}
Sobolev then bounds the positive entropy part; Lemma~\ref{lem:hardy} gives $\int\bigl(\frac{|\nabla f_j|^2}{f_j}+|H_j|^2f_j\bigr)\dd\sigma_j-(m-1)^2b_j^2\ge0$, which bounds the negative part from the upper bound for the functional.

For $s\le1$, ambient homothety gives curvature parameter $sb_j\le B$ and gives the value
$$
s^2\bigl(\int\bigl(\frac{|\nabla f_j|^2}{f_j}+|H_j|^2f_j\bigr)\dd\sigma_j-(m-1)^2b_j^2\bigr)-\int f_j\log f_j\dd\sigma_j-m\log s.
$$
Comparing with $a_B$ at $s=1-\sqrt{\eps_j}$ and expanding yields
\begin{equation}\label{eq:near-scale}
\int\bigl(\frac{|\nabla f_j|^2}{f_j}+|H_j|^2f_j\bigr)\dd\sigma_j-(m-1)^2b_j^2\le m/2+o(1).
\end{equation}
This variation is admissible also when $b_j=B$.

We also need an approximate version of \eqref{eq:pi}, uniform for tests depending on $j$. For $|\phi|+|\nabla\phi|\le C$, evaluate the critical density functional along
$$
\frac{f_j(1+t\phi)^2}{\int f_j(1+t\phi)^2\dd\sigma_j}.
$$
Its first three derivatives on a fixed small interval are bounded by the uniform energy and absolute entropy bounds. Minimality up to $\eps_j$ and the bounded second derivative imply that its first derivative at zero has absolute value at most $C_{m,B,C}\sqrt{\eps_j}$.

Center $\phi$ by subtracting $\int\phi f_j\dd\sigma_j$. Subtract the first derivative along the normalized density perturbation with factor $\bigl[1+t(\phi-\int\phi f_j\dd\sigma_j)^2\bigr]^2$ from the second derivative along the normalized density perturbation with factor $\bigl[1+t(\phi-\int\phi f_j\dd\sigma_j)\bigr]^2$. Direct differentiation gives
$$
8\int|\nabla\phi|^2f_j\dd\sigma_j
-4\Var_{f_jd\sigma_j}(\phi).
$$
The first derivative being subtracted is $O_{m,B,C}(\sqrt{\eps_j})$, since $(\phi-\int\phi f_j\dd\sigma_j)^2$ and its gradient are uniformly bounded. The symmetric Taylor formula at $t=\eps_j^{1/4}$ bounds the second derivative below by $-C_{m,B,C}(\eps_j/t^2+t)=-o(1)$. Hence
\begin{equation}\label{eq:near-pi}
\Var_{f_jd\sigma_j}(\phi)
\le2\int|\nabla\phi|^2f_j+o_{j,m,B,C}(1).
\end{equation}

Suppose a subsequence has $b_j\to b_*>0$, so $\delta\le b_j\le B$ for some $\delta>0$. Use the centers in Lemma~\ref{lem:radial} for $f_jd\sigma_j$. The coordinate functions in \eqref{eq:profile-center} are bounded by $m/((m-1)\delta)$ and have gradients at most one. Apply \eqref{eq:near-pi} to all $m+k$ coordinates and sum. Since $k$ is fixed, their total error is still $o(1)$. This gives \eqref{eq:profile-pi} with an $o(1)$ error, whereas \eqref{eq:profile-div} remains exact. Equations \eqref{eq:radial-algebra} and \eqref{eq:near-scale}, and the bound $\gamma_{b_j}\le m/((m-1)\delta)$, now give
$$
\int\left[b_j^2\gamma_{b_j}^2-
\bigl(b_j\coth(b_jr_j)\gamma_{b_j}-\gamma_{b_j}'\bigr)^2\right]f_j\dd\sigma_j\longrightarrow0.
$$
Indeed the nonnegative right side of \eqref{eq:radial-algebra} is at most
\begin{align*}
&\left(\int\bigl(\frac{|\nabla f_j|^2}{f_j}+|H_j|^2f_j\bigr)\dd\sigma_j-(m-1)^2b_j^2-\frac m2\right)\int\gamma_{b_j}^2f_j\dd\sigma_j\\
&\quad+\frac m2\left\{\int\gamma_{b_j}^2f_j\dd\sigma_j
-2\int\left[(\gamma_{b_j}')^2|\nabla_{\Sigma_j}r_j|^2+
\frac{b_j^2\gamma_{b_j}^2}{\sinh^2(b_jr_j)}(m-|\nabla_{\Sigma_j}r_j|^2)\right]f_j\dd\sigma_j\right\},
\end{align*}
whose upper limit is nonpositive.

Positivity on compact annuli in \eqref{eq:profile-positive} shows that their masses tend to zero. On the other hand, Lemma~\ref{lem:hardy} and \eqref{eq:near-scale} give
$$
\int_{\{r_j<\rho\}}f_j\dd\sigma_j\le\rho^2\bigl(m/2+o(1)\bigr).
$$
First let $j\to\infty$, and then $\rho\downarrow0$. It follows that $\int_{\{r_j\le R\}}f_j\dd\sigma_j\to0$ for every fixed $R$. The uniform limits in Lemma~\ref{lem:profile}, and their bounded integrands, consequently give
\begin{align*}
&\int\gamma_{b_j}^2f_j\dd\sigma_j\longrightarrow\frac{m^2}{(m-1)^2b_*^2}>0,\\
&\int\left[(\gamma_{b_j}')^2|\nabla_{\Sigma_j}r_j|^2+
\frac{b_j^2\gamma_{b_j}^2}{\sinh^2(b_jr_j)}(m-|\nabla_{\Sigma_j}r_j|^2)\right]f_j\dd\sigma_j\longrightarrow0.
\end{align*}
This contradicts the approximate version of \eqref{eq:profile-pi}. Thus $b_j\to0$. Proposition~\ref{prop:unshifted} then implies
$$
a_B=\lim_j\bigl(\int\bigl(\frac{|\nabla f_j|^2}{f_j}+|H_j|^2f_j-f_j\log f_j\bigr)\dd\sigma_j-(m-1)^2b_j^2\bigr)
\ge m+\frac m2\log(4\pi).
$$
Each allowed $b$ belongs to some finite interval $[0,B]$, so the critical inequality holds for all data. For $0\le\beta\le (m-1)^2$, add $((m-1)^2-\beta)b^2\ge0$. Lemma~\ref{lem:density} gives exactly \eqref{eq:main}. The Euclidean Gaussian proves optimality.
\end{proof}

\section{Equality in the logarithmic inequality}\label{sec:equality}

The equality statements concern the geometry carried by the submanifold. A totally geodesic Euclidean submanifold can lie in an ambient manifold with curved directions. For example, a complete geodesic in $\mathbb H^2(-1)$ times $\mathbb R^{m-1}\times\{0\}^{k-1}$ is a totally geodesic Euclidean $m$-submanifold of $\mathbb H^2(-1)\times\mathbb R^{m+k-2}$. The isometry in Theorem~\ref{thm:main} is an isometry of induced length metrics; it gives a smooth totally geodesic parametrization even if the original parametrization is only $C^{1,1}$.

\begin{proof}[Proof of Theorem~\ref{thm:main}: equality for $m\ge3$]
Suppose first that equality holds at $\beta=(m-1)^2$. Since this is the infimum over the entire class, arbitrary ambient homotheties are admissible variations. Their first derivative gives
$$
\int\bigl(\frac{|\nabla f|^2}{f}+|H|^2f\bigr)\dd\sigma-(m-1)^2b^2=m/2.
$$
Normalized multiplicative variations, with the zero-error second variation used in \eqref{eq:near-pi}, give \eqref{eq:pi} with $\kappa=2$. Lemma~\ref{lem:radial} excludes $b>0$. Therefore $b=0$ and $\int\bigl(\frac{|\nabla f|^2}{f}+|H|^2f\bigr)\dd\sigma=m/2$. If $\beta<(m-1)^2$, equality in \eqref{eq:main} first forces equality in the critical inequality and $((m-1)^2-\beta)b^2=0$, and gives the same conclusion.

The source is connected: a component indicator has zero gradient, so \eqref{eq:pi} and positivity of $f$ make every component have mass zero or one. Lemma~\ref{lem:moments} supplies a center $p$. Identify $T_pM$ with $\R^{m+k}$ and consider the Euclidean immersion $\log_p\circ\iota$. The comparison inequalities are
$$
\int r^2f\dd\sigma\le2\int\sum_{\alpha=1}^{m+k}|\nabla (\log_p\circ\iota)_\alpha|^2f\dd\sigma\le2m,
\qquad
m\le\sqrt{\frac m2\int r^2f\dd\sigma}\le m.
$$
All are equalities. In particular $\int r^2f\dd\sigma=2m$, $(\log_p\circ\iota)^*g_{\R^{m+k}}=g_\Sigma$, and equality in Cauchy--Schwarz gives
\begin{equation}\label{eq:equality-vector}
\frac{\nabla f}{f}+H=-\frac{\nabla^M(r^2/2)}2,
\qquad H=-\frac{(\nabla^M(r^2/2))^\perp}{2},
\qquad f=C_0e^{-r^2/4}.
\end{equation}
The metric equality follows pointwise from continuity of $d(\log_p\circ\iota)$ and positivity of the density. The tangential gradient identity gives $\nabla(fe^{r^2/4})=0$ weakly on the connected source, which proves the last formula with $C_0>0$.

We check the mean curvature of the Euclidean immersion $(\log_p\circ\iota)$, since it cannot simply be identified with the ambient mean curvature. Write the ambient metric in global normal coordinates as the positive matrix $G(z)$. Rauch comparison and the Gauss lemma give $G\ge I$ and $G(z)z=z$. The metric equality shows that $G-I$ annihilates $d(\log_p\circ\iota)(T\Sigma)$. For each compactly supported vector field $W$ along $\log_p\circ\iota$ and its variation $\log_p\circ\iota+tW$, its volume density for $G$ is at least its Euclidean volume density, with equality at $t=0$. Their first variations therefore agree. Denote by $H_E$ the mean curvature vector of the Euclidean immersion $\log_p\circ\iota$. We obtain
$$
H_E=G\,d\log_p(H)\quad\hbox{almost everywhere}.
$$
Both the radial vector and all tangent vectors lie in $\ker(G-I)$. Hence $d\log_p((\nabla^M(r^2/2))^\perp)$ is the Euclidean vector $(\log_p\circ\iota)^\perp$ and belongs to this kernel. Equation \eqref{eq:equality-vector} gives
$$
H_E=-(\log_p\circ\iota)^\perp/2.
$$
Thus $(\log_p\circ\iota)$ is a complete Euclidean $C^{1,1}$ self-shrinker, with the same source measure. Normalization and $\int r^2f\dd\sigma=2m$ give
$$
\int\bigl(\frac{|\nabla f|^2}{f}+|H|^2f-f\log f\bigr)\dd\sigma
=m-\log C_0.
$$
Equality therefore implies $C_0=(4\pi)^{-m/2}$ and
$$
(4\pi)^{-m/2}\int e^{-|(\log_p\circ\iota)|^2/4}\dd\sigma=1.
$$

We give the rigidity argument directly in this regularity and immersion class. For $y\in\R^{m+k}$ and $0<\tau<1$, restrict the heat-time parameter to $t=1+(\tau-1)s^2$. Direct differentiation and the weak first variation give
\begin{align}\label{eq:shrinker-monotone}
&\frac{d}{ds}\left[(4\pi t)^{-m/2}
\int e^{-|\log_p\circ\iota-sy|^2/(4t)}\dd\sigma\right]\notag\\
&\quad=-\frac{s}{2t^2}(4\pi t)^{-m/2}
\int|(y+(\tau-1)s(\log_p\circ\iota))^\perp|^2
e^{-|\log_p\circ\iota-sy|^2/(4t)}\dd\sigma\le0.
\end{align}
For the coefficients, multiply the direct derivative by $2t^2$. Its integrand is $-2m(\tau-1)st+\langle (\log_p\circ\iota)-sy,y+(\tau-1)s(\log_p\circ\iota)\rangle$. Subtract the zero weighted divergence identity for the tangent part of $y+(\tau-1)s(\log_p\circ\iota)$, multiplied by $-2t$. The result is $-s|(y+(\tau-1)s(\log_p\circ\iota))^\perp|^2$, using $H_E=-(\log_p\circ\iota)^\perp/2$. This proves \eqref{eq:shrinker-monotone}. For fixed $s>0$, polynomial factors times its Gaussian are bounded by a constant times $e^{-|(\log_p\circ\iota)|^2/4}$ because $t<1$. Complete intrinsic cutoffs therefore justify integration by parts and differentiation. Completing the square bounds the kernel ratio to $e^{-|(\log_p\circ\iota)|^2/4}$ uniformly for $0\le s\le1$, so dominated convergence gives the value one at $s=0$. Hence $(4\pi\tau)^{-m/2}\int e^{-|(\log_p\circ\iota)-y|^2/(4\tau)}\dd\sigma\le1$.

For each source point $x_0$ with $y=\log_p(\iota(x_0))$, one local immersion patch contributes a Gaussian integral tending to one as $\tau\downarrow0$. All other contributions are nonnegative; therefore $\liminf_{\tau\downarrow0}(4\pi\tau)^{-m/2}\int e^{-|(\log_p\circ\iota)-y|^2/(4\tau)}\dd\sigma\ge1$. If $(y-s_0\log_p(\iota(x)))^\perp$ were nonzero at some $x$ and $s_0\in(0,1)$, continuity would give a compact source patch and an interval about $s_0$ on which the integral in \eqref{eq:shrinker-monotone} has a uniform positive lower bound for small $\tau$. Integrating over that interval would give $(4\pi\tau)^{-m/2}\int e^{-|(\log_p\circ\iota)-y|^2/(4\tau)}\dd\sigma\le1-c$ for a fixed $c>0$, a contradiction. Thus
$$
(y-s\log_p(\iota(x)))^\perp=0\quad
\hbox{for every }y\in \log_p(\iota(\Sigma)),\ x\in\Sigma,\ 0<s<1.
$$
Taking two different values of $s$ gives $\log_p(\iota(x))^\perp=0$ and $y^\perp=0$ at every $x$. Fixing $x$ shows that the entire image lies in the linear tangent $m$-plane there. Completeness makes the resulting local isometry onto this plane a covering: straight segments lift for their full length, since finite-length lifts are Cauchy and the source is complete. The plane is simply connected and the source is connected, so this covering has one sheet.

In these coordinates the original immersion is $\exp_p$ restricted to that linear plane, and is smooth. The scalar Gauss equation, its flat induced metric and \eqref{eq:equality-vector} give
$$
0=\sum_{i\ne j}\Sec_M(e_i,e_j)+|H|^2-|II|^2
\le-|II|^2.
$$
Here $H=0$ because $(\log_p\circ\iota)^\perp=0$. Thus $II=0$, and the image is a complete totally geodesic Euclidean $m$-submanifold. Equation \eqref{eq:equality-vector} and normalization give the stated Gaussian. Conversely, direct integration of that Gaussian on any such submanifold gives equality.
\end{proof}

\section{The logarithmic inequality in dimension two}\label{sec:log-two}

Throughout this section $m=2$. Proposition~\ref{prop:unshifted} already gives the unshifted inequality. A product with a Gaussian line also transfers its equality classification from dimension three. The curvature correction then follows from the explicit two-dimensional radial profile and a local divergence estimate.

\begin{proof}[Proof of Theorem~\ref{thm:main} for $m=2$]
Proposition~\ref{prop:unshifted} gives
\begin{equation}\label{eq:two-unshifted}
\int_\Sigma\left(\frac{|\nabla f|^2}{f}+|H|^2f-f\log f\right)\dd\sigma
\ge2+\log(4\pi).
\end{equation}
For equality, use the complete $C^{1,1}$ product immersion $\Sigma\times\R\to M\times\R$, whose mean-curvature vector is $(H,0)$, and the Gaussian $g_1(t)=(4\pi)^{-1/2}e^{-t^2/4}$. The product density has finite energy and absolute entropy. Its Fisher information and entropy are additive, and
$$
\int_\R\frac{|g_1'|^2}{g_1}\dd t=\frac12,\qquad
\int_\R g_1\log g_1\dd t=-\frac{1+\log(4\pi)}2.
$$
Equality in \eqref{eq:two-unshifted} gives equality for the product density in dimension three. Its source is isometric to $\R^3$, and its immersion is a one-sheet totally geodesic embedding. The surface is connected, complete and flat, and $\pi_1(\Sigma)=\pi_1(\Sigma\times\R)=0$, so it is isometric to $\R^2$. The product second fundamental form restricts to $II$ and product injectivity implies injectivity of $\iota$. The natural parallel line field becomes a constant Euclidean direction under the product isometry. Splitting the three-dimensional Gaussian in that direction gives exactly the surface Gaussian in Theorem~\ref{thm:main}. Conversely, that Gaussian attains \eqref{eq:two-unshifted}.

For the variational argument, \eqref{eq:unshifted-homogeneous} with $m=2$ reads
\begin{equation}\label{eq:two-unshifted-hom}
\int_\Sigma f\log f\dd\sigma
\le\log\left\{\frac1{4\pi e}\int_\Sigma\left(\frac{|\nabla f|^2}{f}+|H|^2f\right)\dd\sigma\right\}.
\end{equation}

We will also need a uniform bound on positive entropy. Fix $h\in C_c^1(\R)$ with $\int h^2=1$. Applying Proposition~\ref{thm:sharp-sobolev} in dimension three to $v(x)h(t)$ on the same product gives
\begin{align*}
9|\mathbb S^3|^{2/3}\left(\int_\R |h|^6\dd t\right)^{1/3}\|v\|_6^2
\le{}&\int_\Sigma(12|\nabla v|^2+|H|^2v^2)\dd\sigma\\
&+12\left(\int_\R|h'|^2\dd t\right)\|v\|_2^2.
\end{align*}
Complete intrinsic cutoffs and value truncation extend this inequality to finite-energy $v$. In particular,
$$
\|f\|_3\le C\left[1+\int_\Sigma\left(\frac{|\nabla f|^2}{f}+|H|^2f\right)\dd\sigma\right],
$$
where $C$ depends only on the fixed choice of $h$. This bounds $\int f\log_+f$ when the energy is bounded.

Fix $B<\infty$, and let $a_B$ be the infimum of
$$
\int_\Sigma\left(\frac{|\nabla f|^2}{f}+|H|^2f-f\log f\right)\dd\sigma-b^2
$$
over the stated positive densities and geometries of codimension $k$ with $0\le b\le B$. Inequality~\eqref{eq:two-unshifted-hom} gives $a_B\ge2+\log(4\pi)-B^2$, and the Euclidean Gaussian gives $a_B\le2+\log(4\pi)$. Choose data with values at most $a_B+\eps_j$, where $\eps_j\downarrow0$. Then
\begin{align*}
&\int_{\Sigma_j}\left(\frac{|\nabla f_j|^2}{f_j}+|H_j|^2f_j\right)\dd\sigma_j\\
&\quad-\log\left\{\frac1{4\pi e}\int_{\Sigma_j}\left(\frac{|\nabla f_j|^2}{f_j}+|H_j|^2f_j\right)\dd\sigma_j\right\}
\le3+\log(4\pi)+B^2
\end{align*}
for large $j$. The scalar function $x-\log(x/(4\pi e))$ diverges at both endpoints of $(0,\infty)$, so the energies are uniformly bounded and bounded away from zero. The preceding $L^3$ estimate bounds their positive entropy. The identity expressing their entropy as energy minus $b_j^2$ minus the functional then bounds the negative entropy. Thus the energy and absolute entropy are uniformly bounded in terms of $B$.

For $s\le1$, the homothety remains in the class $b\le B$. Comparing its value with $a_B$ at $s=1-\sqrt{\eps_j}$ gives
\begin{equation}\label{eq:two-near-scale}
\int_{\Sigma_j}\left(\frac{|\nabla f_j|^2}{f_j}+|H_j|^2f_j\right)\dd\sigma_j-b_j^2\le1+o(1).
\end{equation}
Indeed, the rescaled value is the original value plus $(s^2-1)$ times the shifted energy minus $2\log s$; expansion at $s=1$ gives this bound, since the shifted energies are bounded.

The normalized variations used in \eqref{eq:near-pi} also apply in dimension two. For completeness, fix $K>0$ and, for $|\phi|+|\nabla\phi|\le K$, vary $f_j$ through
$$
\frac{f_j(1+t\phi)^2}{\int_{\Sigma_j}f_j(1+t\phi)^2\dd\sigma_j}.
$$
Its first three functional derivatives on a fixed sufficiently small interval are uniformly bounded. The Fisher numerator is quadratic in $t$, with its mixed term controlled by Cauchy--Schwarz; entropy differentiation is justified by $1+t\phi$ being bounded away from zero and the absolute entropy bound. Near minimality implies that the first derivative is $O_{B,K}(\sqrt{\eps_j})$. Center $\phi$ by subtracting its $f_j$-mean. Subtracting the first variation with squared centered test from the second variation with centered test gives exactly
$$
8\int_{\Sigma_j}|\nabla\phi|^2f_j\dd\sigma_j-4\Var_{f_jd\sigma_j}(\phi).
$$
The potential contributions cancel by the weak product rule, and differentiating the normalized entropy gives the variance term. The subtracted first derivative is $O_{B,K}(\sqrt{\eps_j})$. The symmetric Taylor formula at $t=\eps_j^{1/4}$ bounds the second derivative below by $-O_{B,K}(\eps_j^{1/2}+\eps_j^{1/4})$. Consequently,
\begin{equation}\label{eq:two-near-pi}
\Var_{f_jd\sigma_j}(\phi)\le2\int_{\Sigma_j}|\nabla\phi|^2f_j\dd\sigma_j+o_{j,B,K}(1),
\end{equation}
uniformly for these tests, including tests depending on $j$.

Suppose a subsequence has $b_j\to b_*>0$, so $\delta\le b_j\le B$ for some $\delta>0$. The dimension-two version of the radial coefficient in Section~\ref{sec:shift} is
$$
\gamma_b(r)=\frac2b\tanh(br/2),\qquad
\gamma_b'(r)=\operatorname{sech}^2(br/2)=\frac{b\gamma_b(r)}{\sinh(br)}.
$$
We justify its centering without a moment assumption. Fix $p_0\in M$ and minimize
$$
p\longmapsto\int_\Sigma\int_{d_M(p_0,\iota(x))}^{d_M(p,\iota(x))}\gamma_b(t)\dd t\,f(x)\dd\sigma(x).
$$
The inner integral is bounded in absolute value by $2d_M(p,p_0)/b$. Choose a ball about $p_0$ with mass $\alpha>1/2$. Since $\gamma_b(r)\to2/b$, for any sufficiently small $\eta>0$ the contribution of this ball has lower slope $2\alpha(1-\eta)/b$ in $d_M(p,p_0)$; the complementary contribution has slope at least $-2(1-\alpha)/b$. The sum has positive slope. The function is therefore proper and continuous, and Hopf--Rinow gives a minimizer. Differentiation under the integral is justified by the bounded radial derivative and gives
$$
\int_\Sigma\gamma_b(r)\frac{\log_p(\iota(x))}{r}f(x)\dd\sigma(x)=0.
$$
The integrand has its smooth extension at the pole.

Use these centers $p_j$, and write $r_j=d_M(p_j,\iota_j(\cdot))$. Rauch comparison bounds both the radial and transverse singular values of the differential of $\gamma_b(r)\log_p/r$ by $\operatorname{sech}^2(br/2)$. Its $2+k$ coordinates have absolute value at most $2/\delta$ and gradients at most one. Summing \eqref{eq:two-near-pi} over these centered coordinates, with $k$ fixed, gives
\begin{equation}\label{eq:two-radial-pi}
\int_{\Sigma_j}\gamma_{b_j}(r_j)^2f_j\dd\sigma_j
\le4\int_{\Sigma_j}\operatorname{sech}^4(b_jr_j/2)f_j\dd\sigma_j+o(1).
\end{equation}
For the bounded smooth ambient field $\gamma_b(r)\nabla^Mr$, Hessian comparison gives
$$
\tr_{T\Sigma}\nabla^M(\gamma_b\nabla^Mr)
\ge2+2\tanh^2(br/2)(1-|\nabla_\Sigma r|^2)\ge2.
$$
Intrinsic cutoffs and Fatou first establish integrability of its nonnegative trace. Equation~\eqref{eq:weighted-div} and Cauchy--Schwarz then give the exact estimate
\begin{equation}\label{eq:two-radial-div}
4\le\int_{\Sigma_j}\left(\frac{|\nabla f_j|^2}{f_j}+|H_j|^2f_j\right)\dd\sigma_j
\int_{\Sigma_j}\gamma_{b_j}(r_j)^2f_j\dd\sigma_j.
\end{equation}
Combining \eqref{eq:two-near-scale}--\eqref{eq:two-radial-div}, using $\gamma_{b_j}\le2/\delta$ and $b_j^2\gamma_{b_j}^2=4\tanh^2(b_jr_j/2)$, yields
\begin{align*}
4&\le b_j^2\int_{\Sigma_j}\gamma_{b_j}^2f_j\dd\sigma_j
+4\int_{\Sigma_j}\operatorname{sech}^4(b_jr_j/2)f_j\dd\sigma_j+o(1)\\
&=4-4\int_{\Sigma_j}\tanh^2(b_jr_j/2)\operatorname{sech}^2(b_jr_j/2)f_j\dd\sigma_j+o(1).
\end{align*}
The last integrand has a uniform positive lower bound on each fixed annulus $0<\rho\le r_j\le R$. Thus every such annulus has mass tending to zero.

To control the mass near the pole, choose a smooth nonincreasing function $\chi:[0,\infty)\to[0,1]$, equal to one on $[0,1]$, zero on $[2,\infty)$, with $|\chi'|\le2$. Hessian comparison applied to $Y=\chi(r/\rho)\nabla^M(r^2/2)$ gives
$$
\tr_{T\Sigma}\nabla^MY\ge2\mathbf1_{\{r\le\rho\}}-4\mathbf1_{\{\rho<r<2\rho\}},
\qquad |Y|\le2\rho.
$$
The field is smooth at the pole. Complete intrinsic cutoffs make its boundary term tend to zero, also for a nonproper immersion. Equation~\eqref{eq:weighted-div} consequently implies
$$
\int_{\{r\le\rho\}}f\dd\sigma
\le\rho\left[\int_\Sigma\left(\frac{|\nabla f|^2}{f}+|H|^2f\right)\dd\sigma\right]^{1/2}
+2\int_{\{\rho<r<2\rho\}}f\dd\sigma.
$$
Apply this to the sequence. For each fixed $R$, first let $j\to\infty$ and then $\rho\downarrow0$; the annulus conclusion and the uniform energy bound show that the mass of $\{r_j\le R\}$ tends to zero. Uniform radial limits for $\delta\le b_j\le B$ now give
$$
\int_{\Sigma_j}\gamma_{b_j}^2f_j\dd\sigma_j\longrightarrow\frac4{b_*^2},
\qquad
\int_{\Sigma_j}\operatorname{sech}^4(b_jr_j/2)f_j\dd\sigma_j\longrightarrow0.
$$
This contradicts \eqref{eq:two-radial-pi}. Therefore $b_j\to0$, and \eqref{eq:two-unshifted} gives $a_B\ge2+\log(4\pi)$. As $B$ was arbitrary, \eqref{eq:main} holds for $m=2$, $\beta=1$. Adding $(1-\beta)b^2$ proves it for $0\le\beta\le1$.

If equality holds at $\beta=1$, arbitrary homotheties in both directions give shifted energy exactly one. The same normalized variations give the exact Poincar\'e inequality with constant two. If $b>0$, the exact versions of \eqref{eq:two-radial-pi} and \eqref{eq:two-radial-div} at the preceding center imply
$$
0\ge4\int_\Sigma\tanh^2(br/2)\operatorname{sech}^2(br/2)f\dd\sigma>0.
$$
The strict inequality follows because $f>0$ almost everywhere and the pole fiber is discrete and has zero source area. Thus $b=0$, and the equality classification following \eqref{eq:two-unshifted} applies. Equality for $\beta<1$ first forces equality for $\beta=1$ and gives the same conclusion. The Euclidean Gaussian establishes optimality of the additive constant.
\end{proof}

\section{The logarithmic threshold and homogeneous form}\label{sec:sharpness}

In this section $m\ge2$. The following proposition completes the sharpness assertion for the curvature coefficient in Theorem~\ref{thm:main}.

\begin{proposition}\label{prop:threshold}
For every $\beta>(m-1)^2$ and fixed $k$, the left side of \eqref{eq:main} has no finite lower bound depending only on $m$ and $\beta$, uniformly over $b\ge0$ and all data satisfying the geometric and density hypotheses of Theorem~\ref{thm:main}. For each fixed $b\ge0$, all such data satisfy
\begin{equation}\label{eq:supercritical-bound}
\int_\Sigma\left(\frac{|\nabla f|^2}{f}-f\log f+
(|H|^2-\beta b^2)f\right)\dd\sigma
\ge m+\frac m2\log(4\pi)-[\beta-(m-1)^2]b^2.
\end{equation}
The estimate \eqref{eq:supercritical-bound} does not determine the best lower bound for each fixed $b>0$.
\end{proposition}

\begin{proof}[Proof of Proposition~\ref{prop:threshold}]
By normalization, the left side of \eqref{eq:supercritical-bound} equals the left side of \eqref{eq:main} at $\beta=(m-1)^2$, minus $[\beta-(m-1)^2]b^2$. Theorem~\ref{thm:main} therefore proves \eqref{eq:supercritical-bound} as well. For $\beta=m^2$, its right side is $m+(m/2)\log(4\pi)-(2m-1)b^2$.

To exclude a lower bound independent of $b$, realize $\mathbb H^m(-b^2)$ as a totally geodesic submanifold of $\mathbb H^{m+k}(-b^2)$. Its trace mean curvature vanishes. Choose $a>(m-1)$ close enough to $(m-1)$ that $ma^2/(a+1)<\beta$. On this submanifold take
$$
f_{b,a}(r)=\frac{\Gamma((a+1)/2)}{\pi^{m/2}\Gamma((a-m+1)/2)}
b^m\cosh(br)^{-a}.
$$
The polar volume element is $(\sinh(br)/b)^{m-1}\dd r\dd\omega$. The substitution $t=\tanh(br)$ reduces normalization to the beta integral
$$
\frac12\int_0^1 v^{m/2-1}(1-v)^{(a-m-1)/2}\dd v.
$$
It also gives $\int\tanh^2(br)f_{b,a}\dd\sigma=m/(a+1)$. Consequently
$$
\int\frac{|\nabla f_{b,a}|^2}{f_{b,a}}\dd\sigma
=\frac{ma^2}{a+1}b^2,
\qquad
\int f_{b,a}\log f_{b,a}\dd\sigma=m\log b+\int_{\mathbb H^m(-1)}f_{1,a}\log f_{1,a}\dd\sigma_1.
$$
Here $d\sigma_1$ is the volume measure at curvature $-1$. The exponential tail with $a>m-1$ gives finite absolute entropy and Fisher information, so the two integral conditions hold. The left side of \eqref{eq:main} equals
$$
\left(\frac{ma^2}{a+1}-\beta\right)b^2-m\log b-\int_{\mathbb H^m(-1)}f_{1,a}\log f_{1,a}\dd\sigma_1,
$$
which tends to $-\infty$ as $b\to\infty$.
\end{proof}

\begin{corollary}\label{cor:homogeneous}
Under the hypotheses of Theorem~\ref{thm:main}, for every $0\le\beta\le(m-1)^2$,
\begin{equation}\label{eq:homogeneous}
\int_\Sigma f\log f\dd\sigma
\le\frac m2\log\left\{
\frac1{2\pi e m}\int_\Sigma\left[
\frac{|\nabla f|^2}{f}+(|H|^2-\beta b^2)f\right]\dd\sigma\right\}.
\end{equation}
The integral inside the logarithm is strictly positive. Equality holds precisely on the same totally geodesic Euclidean submanifolds with $b=0$, for the densities
$$
f\circ\Phi^{-1}(x)=(4\pi\tau)^{-m/2}e^{-|x-x_0|^2/(4\tau)},
\qquad \tau>0.
$$
\end{corollary}
\begin{proof}
Apply Theorem~\ref{thm:main} to the ambient metric $s^{-2}g$, curvature parameter $sb$ and probability density $s^mf$. Its conclusion is
$$
s^2\int_\Sigma\left[\frac{|\nabla f|^2}{f}
+(|H|^2-\beta b^2)f\right]\dd\sigma
-\int_\Sigma f\log f\dd\sigma-m\log s
\ge m+\frac m2\log(4\pi),\qquad s>0.
$$
If the shifted energy integral were nonpositive, letting $s\to\infty$ would contradict this bound. Thus that integral is strictly positive. The minimum on the left occurs at
$$
s^2=\frac{m}{2\int_\Sigma\left[
\frac{|\nabla f|^2}{f}+(|H|^2-\beta b^2)f\right]\dd\sigma},
$$
and gives \eqref{eq:homogeneous}. Equality is equivalent to equality in the additive inequality on this rescaled geometry. Its classification gives $b=0$ and the stated Euclidean submanifold. Returning to the original metric gives a Gaussian with $\tau=s^2>0$. Every such density satisfies
$$
\int_\Sigma\left[\frac{|\nabla f|^2}{f}
+(|H|^2-\beta b^2)f\right]\dd\sigma=\frac m{2\tau},
$$
and attains \eqref{eq:homogeneous}.
\end{proof}

We record the extension to densities that vanish, needed for the Nash inequality. It also applies to the higher-dimensional logarithmic inequality.

\begin{lemma}\label{lem:nonnegative}
For $m\ge2$, Theorem~\ref{thm:main} and Corollary~\ref{cor:homogeneous} remain valid as inequalities for nonnegative probability densities of finite energy and finite absolute entropy, with the conventions of Lemma~\ref{lem:density}.
\end{lemma}
\begin{proof}
Every source in the stated class carries a positive admissible probability density. To construct one, choose a countable locally finite collection of nonnegative compactly supported smooth functions $h_l$, with $\int_\Sigma h_l^2\dd\sigma=1$, whose positive sets cover the source. Their squares have finite Fisher information, mean-curvature energy and absolute entropy. Choose positive masses $p_l$ with sum one such that
$$
\sum_l p_l\left[\int_\Sigma\left(4|\nabla h_l|^2+|H|^2h_l^2+h_l^2|\log h_l^2|\right)\dd\sigma+|\log p_l|\right]<\infty.
$$
Such masses result by normalizing $2^{-l}$ divided by one plus the displayed integral; the terms involving logarithms are summable since $\log x/x$ is bounded for $x\ge1$. The locally finite sum $h=\sum_l p_lh_l^2$ is positive and smooth. Convexity of Fisher information and linearity of the curvature term give finite energy. Convexity of $t\log_+t$ bounds its positive entropy. Pointwise,
$$
h\log h\ge\sum_l p_lh_l^2\log(p_lh_l^2),
$$
and the right side is absolutely integrable by the choice of the masses. Hence its negative entropy is finite as well.

For any nonnegative admissible $f$, the mixtures $(1-\eps)f+\eps h$ are positive admissible densities. Fisher convexity, linearity of the curvature term, and the bound
\begin{align*}
&\int((1-\eps)f+\eps h)\log((1-\eps)f+\eps h)\dd\sigma\\
&\quad\ge(1-\eps)\int f\log f\dd\sigma+\eps\int h\log h\dd\sigma\\
&\qquad+(1-\eps)\log(1-\eps)+\eps\log\eps
\end{align*}
show that the upper limit of their additive functionals is at most the functional of $f$. This proves the additive inequality for $f$. Applying it after every homothety, as in the proof of Corollary~\ref{cor:homogeneous}, gives strict positivity of its shifted energy and the homogeneous inequality. For a disconnected source, the bump family covers all components; local finiteness and the summability above justify the same argument.
\end{proof}

\section{Critical exponential inequalities on immersed surfaces}\label{sec:critical-two}

Throughout this section $m=2$. A single explicit radial supersolution gives the heat-kernel bound for both the Dirichlet and global realizations of $-4\Delta+|H|^2-b^2$. We then prove two abstract consequences of this bound: a finite-measure estimate and a global estimate with a positive mass term. The normal-square calculation is related to the entropy-kernel calculation in \cite[Section~5]{BB25}; here the comparison is on the immersed source and retains its multiplicities.

\subsection{Geometric heat-kernel comparison and the global form}

\begin{proposition}\label{global:prop:heat}
The Friedrichs realization
$$
A=-4\Delta_\Sigma+|H|^2-b^2
$$
is nonnegative, and its quadratic form domain is exactly
$$
\mathcal D(A^{1/2})
=\left\{v\in W^{1,2}(\Sigma):
\int_\Sigma |H|^2|v|^2\dd\sigma<\infty\right\}.
$$
For each relatively compact smooth source domain $U$, its Dirichlet realization $A_U$ has compact resolvent and strictly positive lowest eigenvalue. Whole compact source components may be included, with no boundary condition on those components. The global operator has a nonnegative symmetric heat kernel satisfying
\begin{equation}\label{global:eq:offdiag}
K_A(t,x,y)\le
\frac1{16\pi t}
\left(\frac{b\,d_M(\iota(x),\iota(y))}
{\sinh\!\bigl(b\,d_M(\iota(x),\iota(y))\bigr)}\right)^{1/2}
e^{-d_M(\iota(x),\iota(y))^2/(16t)}.
\end{equation}
The quotient in parentheses is interpreted as $1$ if its argument is zero. Its canonical diagonal obeys
\begin{equation}\label{global:eq:diag}
K_A(t,x,x)=
\int_\Sigma |K_A(t/2,x,z)|^2\dd\sigma(z)
\le\frac1{16\pi t},\qquad t>0.
\end{equation}
\end{proposition}

\subsubsection{A radial supersolution}

Define the comparison function
$$
G_b(t,r)=\frac1{16\pi t}
\left(\frac{br}{\sinh(br)}\right)^{1/2}e^{-r^2/(16t)},
\qquad t>0,\quad r\ge0,
$$
with the same limiting convention. For $b>0$ and $r>0$,
\begin{equation}\label{global:eq:logderivative}
(\log G_b)_r=\frac1{2r}-\frac b2\coth(br)-\frac r{8t}<0.
\end{equation}
Direct differentiation gives
\begin{align}
&\partial_t\log G_b
-4\left[(\log G_b)_{rr}
+b\coth(br)(\log G_b)_r+((\log G_b)_r)^2\right]-b^2
\notag\\
&\hspace{35mm}=\frac1{r^2}-\frac{b^2}{\sinh^2(br)}\ge0.
\label{global:eq:scalarresidual}
\end{align}
For completeness, the time-independent part of $\log G_b$ has derivative
$(r^{-1}-b\coth(br))/2$. Thus the terms of order $t^{-1}$ cancel, and the remaining identity follows by expanding this derivative and its square.

The second scalar inequality needed below is
\begin{align}
(\log G_b)_{rr}-b\coth(br)(\log G_b)_r
={}&\frac{br\coth(br)-1}{8t}\notag\\
&+\frac{b^2}{2}
\left(1+\frac2{\sinh^2(br)}
-\frac{\coth(br)}{br}-\frac1{(br)^2}\right)\ge0.
\label{global:eq:radialconvex}
\end{align}
To verify the last sign, for every $x>0$,
\begin{align*}
&2x^2\sinh^2x
\left(1+\frac2{\sinh^2x}-\frac{\coth x}{x}-\frac1{x^2}\right)\\
&\quad=(x^2-1)\cosh(2x)-x\sinh(2x)+3x^2+1\\
&\quad=\sum_{j=3}^{\infty}
\frac{2^{2j-2}(2j+1)(2j-4)}{(2j)!}x^{2j}\ge0.
\end{align*}
The power series converges everywhere and all displayed coefficients are positive. Also $x\coth x\ge1$. These facts prove \eqref{global:eq:radialconvex}. When $b=0$, replace $b\coth(br)$ by $1/r$; both residuals in \eqref{global:eq:scalarresidual} and \eqref{global:eq:radialconvex} vanish.

Fix $o\in M$, and put $r=d_M(\iota(\cdot),o)$. The Cartan--Hadamard theorem removes the cut locus of $o$. Hessian comparison under the upper sectional-curvature bound gives, off the fiber of $o$,
\begin{equation}\label{global:eq:hessian}
\tr_{T\Sigma}\Hess_M r
\ge b\coth(br)\left(2-|\nabla_\Sigma r|^2\right).
\end{equation}
The comparison has the required direction: the index form of a perpendicular radial Jacobi field is bounded below by
$\int_0^r(|D_sJ|^2+b^2|J|^2)\dd s$. Minimizing this latter integral with prescribed endpoint gives $b\coth(br)$ times the squared endpoint length. Taking the trace over the source tangent plane gives \eqref{global:eq:hessian}.

The immersion chain rule holds weakly:
$$
\Delta_\Sigma(\psi\circ\iota)
=\tr_{T\Sigma}\Hess_M\psi+
\langle H,\nabla^M\psi\rangle.
$$
Since \eqref{global:eq:logderivative} is negative, \eqref{global:eq:hessian} gives the upper bound needed when the Laplacian of $G_b(t,r)$ is expanded. Completing the full normal-vector square yields
\begin{align}
\frac{(\partial_t-4\Delta_\Sigma+|H|^2-b^2)G_b(t,r)}
{G_b(t,r)}
\ge{}&
\frac1{r^2}-\frac{b^2}{\sinh^2(br)}
\notag\\
&+\left|H-2(\log G_b)_r(\nabla^Mr)^\perp\right|^2
\notag\\
&+4(1-|\nabla_\Sigma r|^2)
\left[(\log G_b)_{rr}
-b\coth(br)(\log G_b)_r\right]
\ge0.
\label{global:eq:supersolution}
\end{align}
Here $|(\nabla^Mr)^\perp|^2=1-|\nabla_\Sigma r|^2$. In particular, the square is the squared norm of a normal vector in the actual codimension; no estimate of its components introduces a factor depending on $k$.

For each $t>0$, $G_b(t,d_M(\cdot,o))$ is smooth at $o$: squared distance is smooth there and the radial expression is a smooth even function. Its composition with a $C^{1,1}$ immersion is locally $C^{1,1}$. Each fiber is locally discrete, since an immersion is locally an embedding. The induced metric and its volume density are locally Lipschitz, so the weak divergence-form Laplacian of the composite has a locally bounded density. Consequently \eqref{global:eq:supersolution} holds weakly across all preimages of $o$, with no additional point masses.

\subsubsection{Comparison on compact source domains}

Let $U\Subset\Sigma$ be a smooth source domain; whole compact components may also be included. On its compact closure the induced metric is locally uniformly elliptic and Lipschitz, and $|H|^2-b^2$ is bounded. The Dirichlet form on $W^{1,2}_0(U)$ is closed and bounded below by $-b^2\|\cdot\|_2^2$; for a whole compact component the form domain is $W^{1,2}(U)$. The associated operator $A_U$ has compact resolvent.

The local heat-kernel theory applies to these coefficients. Norris \cite[pp.~79--81]{No97} allows measurable locally uniformly elliptic metric coefficients and a locally comparable measure, with no global uniform ellipticity assumption. Our coefficients satisfy these hypotheses. A bounded real potential is included by the bounded-perturbation expansion; positivity follows by the negative-part maximum principle, equivalently by positive Trotter factors. Interior parabolic regularity supplies a continuous kernel for positive times. Symmetry and the semigroup law give its canonical diagonal.

Fix $y\in U$ and set $o=\iota(y)$. In source coordinates centered at $y$ and orthonormal at that point, the $C^{1,1}$ expansion gives
$$
r(z)^2=|z|^2+O(|z|^3),\qquad
d\sigma(z)=(1+O(|z|))\dd z.
$$
Choose a small coordinate neighborhood $N\Subset U$ on which $\iota$ is an embedding. On it, $r(z)\ge c|z|$ for a positive local constant $c$. For every continuous test $\eta$ supported there, the substitution $z=\sqrt\varepsilon\,w$ and Gaussian domination imply
\begin{equation}\label{global:eq:initialmass}
\lim_{\varepsilon\downarrow0}
\int_N G_b(\varepsilon,r(z))\eta(z)\dd\sigma(z)=\eta(y),
\end{equation}
because
$$
\int_{\R^2}\frac1{16\pi}e^{-|w|^2/16}\dd w=1.
$$
The factor $(br/\sinh(br))^{1/2}$ tends to $1$ on the diffusive scale and is bounded by $1$. The distance is bounded away from zero on compact subsets of $N\setminus\{y\}$, so no mass is lost in \eqref{global:eq:initialmass}.

For $0<\varepsilon<t$, the supersolution inequality and the weak maximum principle give
\begin{equation}\label{global:eq:epscomparison}
G_b(t,r(x))
\ge\int_U K_{A_U}(t-\varepsilon,x,z)
G_b(\varepsilon,r(z))\dd\sigma(z).
\end{equation}
The difference has zero initial data at time $\varepsilon$ and nonnegative boundary values. Multiplication of the difference by $e^{-b^2(t-\varepsilon)}$ changes its potential to $|H|^2$, after which testing its negative part proves comparison. All quantities are bounded on compact positive time intervals.

Restrict the right side of \eqref{global:eq:epscomparison} to the selected neighborhood $N$ and use \eqref{global:eq:initialmass} and positive-time continuity. Letting $\varepsilon\downarrow0$ gives
\begin{equation}\label{global:eq:compactkernel}
0\le K_{A_U}(t,x,y)\le G_b(t,d_M(\iota(x),\iota(y))).
\end{equation}
This argument requires unit mass only on the selected source sheet. Other sheets contribute nonnegative quantities in \eqref{global:eq:epscomparison}; their total initial mass need not equal one. Thus self-intersections and nonproperness do not affect the comparison.

In particular, the semigroup identity and \eqref{global:eq:compactkernel} give
$$
e^{-t\lambda_1(A_U)}
\le\tr(e^{-tA_U})
=\int_U K_{A_U}(t,x,x)\dd\sigma(x)
\le\frac{|U|}{16\pi t}.
$$
The trace identity follows from the Hilbert--Schmidt identity for $e^{-tA_U/2}$. The bound makes the required integral finite. Letting $t\to\infty$ proves $\lambda_1(A_U)>0$. Hence
\begin{equation}\label{global:eq:nonnegativecore}
\int_\Sigma\left[4|\nabla v|^2+(|H|^2-b^2)|v|^2\right]\dd\sigma\ge0
\qquad(v\in C_c^\infty(\Sigma)).
\end{equation}

\subsubsection{The exact global form and exhaustion}

The sum of the nonnegative gradient form and the multiplication form,
$$
\int_\Sigma(4|\nabla v|^2+|H|^2|v|^2)\dd\sigma,
$$
is closed on $\{v\in W^{1,2}(\Sigma):Hv\in L^2(\Sigma)\}$. Indeed, a Cauchy sequence for this form plus the $L^2$ norm converges in $W^{1,2}$, and its products with $H$ converge in $L^2$; almost-everywhere subsequences identify the latter limit with the product of $H$ and the former limit. Subtracting the bounded form $b^2\|v\|_2^2$ preserves closedness and lower semiboundedness.

This natural domain is the closure of $C_c^\infty(\Sigma)$ in the norm
$$
\left(\int_\Sigma
[4|\nabla v|^2+(1+|H|^2)|v|^2]\dd\sigma\right)^{1/2}.
$$
To see this, on each connected component use intrinsic distance cutoffs $\chi_R$ with compact support, $0\le\chi_R\le1$, $\chi_R\to1$, and $|\nabla\chi_R|\le c/R$. The source is a complete locally compact length space; the length-space Hopf--Rinow theorem makes closed intrinsic balls compact. The discarded $L^2$, gradient and $H$-weighted tails tend to zero. The additional gradient term is bounded in $L^2$ by $cR^{-1}\|v\|_2$. Each compactly supported Sobolev function can then be approximated in finitely many source charts by smooth functions, because the metric is locally uniformly elliptic and $H$ is bounded on the compact support. On a disconnected source there are countably many components; first retain finitely many components with arbitrarily small omitted form-norm tail, and then use the same construction on each retained component.

This proves the core assertion. Equation~\eqref{global:eq:nonnegativecore} therefore extends to the whole form domain and defines a nonnegative Friedrichs operator $A$ with exactly the domain stated in Proposition~\ref{global:prop:heat}.

Choose smooth relatively compact source domains $U_j\nearrow\Sigma$, retaining whole compact components when needed. Extend their heat kernels by zero outside $U_j\times U_j$. The negative-part maximum principle gives domain monotonicity, so these kernels increase with $j$. The compact comparison \eqref{global:eq:compactkernel} gives a finite pointwise limit.

We check that this limit is the heat kernel of the specified global operator. Extend the Dirichlet resolvent $(A_{U_j}+1)^{-1}$ by restricting its input to $U_j$ and extending its output by zero. Its variational minimization problem is the global coercive form restricted to the closed subspace $W^{1,2}_0(U_j)$. The increasing union of these subspaces is form-dense by the core assertion. Therefore these extended resolvents converge strongly to $(A+1)^{-1}$. This is a statement about extended resolvents, not about a densely defined zero-extended Dirichlet generator on $L^2(\Sigma)$.

The corresponding extended heat operators converge strongly as well. In fact, on the resolvent spectrum $[0,1]$ use the continuous function which is zero at $s=0$ and equals $e^{-t(s^{-1}-1)}$ for $s>0$, and approximate it uniformly by polynomials. Strong convergence of bounded resolvents then gives strong convergence of their heat operators. For a nonnegative compactly supported $L^2$ input, monotone convergence of the kernel integrals identifies the limiting kernel with $e^{-tA}$.

Finally the semigroup identity passes to the limit by monotone convergence. In particular,
\begin{align*}
K_A(t,x,x)
&=\lim_{j\to\infty}K_{A_{U_j}}(t,x,x)\\
&=\lim_{j\to\infty}
\int_{U_j}|K_{A_{U_j}}(t/2,x,z)|^2\dd\sigma(z)\\
&=\int_\Sigma |K_A(t/2,x,z)|^2\dd\sigma(z)
\le\frac1{16\pi t}.
\end{align*}
The full off-diagonal comparison also passes to the limit. This completes the proof of Proposition~\ref{global:prop:heat}. The global resolvent need not be compact, and no global heat trace has been used.

\subsection{Abstract critical exponential estimates}

The two spectral arguments below use the same consequence of Adams's integral lemma.

\begin{lemma}\label{two:lem:oned}
Let $d\ge0$, let $E_0\ge0$, and let $\vartheta:[0,\infty)\to[0,\infty)$ be measurable with $E_0+\int_0^\infty\vartheta^2\le1$. There is a finite constant $B_d$, depending only on $d$, such that
$$
\int_0^\infty
\exp\!\left[\left(\int_0^t\vartheta(s)\dd s
+d\sqrt{E_0+\int_t^\infty\vartheta(s)^2\dd s}\right)^2-t\right]\dd t
\le B_d.
$$
\end{lemma}

\begin{proof}
Apply Adams's integral lemma in the form stated in \cite[Lemma~4.1]{LY16}, with exponent two. Extend $\vartheta$ to $(-\infty,0)$ by setting it equal to $\sqrt{E_0}$ on $(-1,0)$ and zero on $(-\infty,-1]$. For each $t\ge0$, take the kernel equal to one on $(0,t)$ and equal to
$$
\frac{d\vartheta(s)}{\sqrt{E_0+\int_t^\infty\vartheta^2}}
$$
on its complement, using zero when the denominator vanishes. The squared kernel integral over that complement is at most $d^2$, and its pairing with $\vartheta$ is exactly $d\sqrt{E_0+\int_t^\infty\vartheta^2}$. The cited lemma has a constant depending only on this kernel bound and the exponent. It therefore gives the assertion uniformly in $E_0$ and $\vartheta$.
\end{proof}

\begin{proposition}\label{two:prop:abstract}
Let $(X,\mu)$ have finite positive measure, and let $L$ be a lower semibounded self-adjoint operator on $L^2(X,\mu)$ with compact resolvent. Suppose its heat semigroup has a symmetric kernel whose canonical diagonal satisfies
$$
K_L(t,x,x)=\int_X|K_L(t/2,x,y)|^2\dd\mu(y)\le\frac at
\qquad(t>0)
$$
for almost every $x$, where $a>0$. Then $L$ is strictly positive, and there is an absolute finite constant $\mathcal C$ such that every $u$ in its quadratic form domain satisfies
$$
\|L^{1/2}u\|_2^2\le1
\quad\Longrightarrow\quad
\int_X e^{|u|^2/a}\dd\mu\le\mathcal C\mu(X).
$$
The constant $\mathcal C$ is independent of $X,\mu,L$ and $a$.
\end{proposition}

\begin{proof}
The diagonal identity makes $e^{-tL/2}$ Hilbert--Schmidt and gives
$$
\tr(e^{-tL})=\int_XK_L(t,x,x)\dd\mu(x)\le a\mu(X)/t.
$$
Hence the lowest eigenvalue $\lambda_1$ is positive and
\begin{equation}\label{two:eq:gap}
\lambda_1a\mu(X)\ge1/e.
\end{equation}
Indeed, $e^{-t\lambda_1}\le a\mu(X)/t$ first excludes $\lambda_1\le0$; then take $t=1/\lambda_1$. The spectral theorem also gives
\begin{equation}\label{two:eq:largetime}
K_L(t,x,x)\le e^{-\lambda_1t/2}K_L(t/2,x,x)
\le\frac{2a}{t}e^{-\lambda_1t/2}.
\end{equation}
To obtain the first inequality, factor the heat operator and apply its $L^2$ spectral norm to its kernel row.

First take a finite spectral sum $u$ with $\|L^{1/2}u\|_2^2\le1$. Spectral calculus gives
\begin{equation}\label{two:eq:frame}
\int_0^\infty\|Le^{-tL/2}u\|_2^2\dd t=\|L^{1/2}u\|_2^2\le1,
\qquad
u=\int_0^\infty Le^{-tL}u\dd t.
\end{equation}
For an eigenvalue $\lambda>0$, the scalar integrals are $\int_0^\infty\lambda^2e^{-t\lambda}\dd t=\lambda$ and $\int_0^\infty\lambda e^{-t\lambda}\dd t=1$. For a split time $\tau>0$, the portions of the second integral over $(0,\tau)$ and $(\tau,\infty)$ equal $(1-e^{-\tau L})u$ and $e^{-\tau L}u$. Since $(1-e^{-x})^2\le x(1-e^{-x})$,
\begin{equation}\label{two:eq:high}
\|(1-e^{-\tau L})u\|_2^2
\le\tau\int_0^\tau\|Le^{-tL/2}u\|_2^2\dd t.
\end{equation}
For the complementary portion, write $Le^{-tL}u=e^{-tL/2}(Le^{-tL/2}u)$ and use Cauchy--Schwarz against the kernel row to obtain
\begin{equation}\label{two:eq:low}
|e^{-\tau L}u(x)|
\le\int_\tau^\infty\sqrt{K_L(t,x,x)}\,\|Le^{-tL/2}u\|_2\dd t.
\end{equation}

Let $u^*$ be the nonincreasing rearrangement of $|u|$ on $(0,\mu(X))$. For $T\ge0$, evaluate the rearrangement at $\mu(X)e^{-T}$ and take $\tau=a\mu(X)e^{-T}$. Chebyshev's inequality bounds the rearrangement of a bounded function plus an $L^2$ function at $\mu(X)e^{-T}$ by the first function's $L^\infty$ norm plus $e^{T/2}\mu(X)^{-1/2}$ times the second function's $L^2$ norm. Apply this observation to \eqref{two:eq:high} and \eqref{two:eq:low}. Use the original diagonal estimate for $\tau\le t\le a\mu(X)$ and \eqref{two:eq:largetime} for $t\ge a\mu(X)$. Substitution $t=a\mu(X)e^{-S}$, together with \eqref{two:eq:gap}, gives
\begin{align}\label{two:eq:rearrange}
\frac{u^*(\mu(X)e^{-T})}{\sqrt a}
\le{}&\int_{-\infty}^0\sqrt2\,e^{-e^{-S}/(4e)}
\sqrt{a\mu(X)e^{-S}}\,\|Le^{-a\mu(X)e^{-S}L/2}u\|_2\dd S\notag\\
&+\int_0^T\sqrt{a\mu(X)e^{-S}}\,\|Le^{-a\mu(X)e^{-S}L/2}u\|_2\dd S\notag\\
&+\left[\int_T^\infty a\mu(X)e^{-S}\|Le^{-a\mu(X)e^{-S}L/2}u\|_2^2\dd S\right]^{1/2}.
\end{align}
The last coefficient is one because $\tau=a\mu(X)e^{-T}$.

The squared $L^2(-\infty,0)$ norm of the weight $\sqrt2e^{-e^{-S}/(4e)}$ is bounded by
$$
2\int_1^\infty e^{-r/(2e)}\frac{\dd r}{r}\le4e.
$$
Cauchy--Schwarz on this half-line and then in two components bounds the first and last terms in \eqref{two:eq:rearrange} by
$$
\sqrt{4e+1}\left[
\int_{(-\infty,0)\cup(T,\infty)}a\mu(X)e^{-S}\|Le^{-a\mu(X)e^{-S}L/2}u\|_2^2\dd S
\right]^{1/2}.
$$
In Lemma~\ref{two:lem:oned}, substitute
\begin{align*}
\vartheta(S)&=\sqrt{a\mu(X)e^{-S}}\,\|Le^{-a\mu(X)e^{-S}L/2}u\|_2,
\qquad S\ge0,\\
E_0&=\int_{-\infty}^0a\mu(X)e^{-S}\|Le^{-a\mu(X)e^{-S}L/2}u\|_2^2\dd S,
\qquad d=\sqrt{4e+1}.
\end{align*}
These substitutions satisfy the lemma's normalization by \eqref{two:eq:frame}. Equimeasurability and $s=\mu(X)e^{-T}$ now give
\begin{align*}
\int_Xe^{|u|^2/a}\dd\mu
&=\mu(X)\int_0^\infty
\exp\!\left[\frac{u^*(\mu(X)e^{-T})^2}{a}-T\right]\dd T\\
&\le B_{\sqrt{4e+1}}\mu(X).
\end{align*}
For a general form-domain $u$, finite spectral truncations have energies at most one and converge in $L^2$. An almost-everywhere convergent subsequence and Fatou's lemma prove the same conclusion.
\end{proof}

\begin{proposition}\label{global:prop:abstract}
Let $(X,\mu)$ be a $\sigma$-finite measure space, and let $A\ge0$ be self-adjoint on $L^2(X,\mu)$. Suppose its heat operator has a measurable symmetric kernel with canonical diagonal
$$
K_A(t,x,x)=\int_X|K_A(t/2,x,y)|^2\dd\mu(y)\le\frac at
$$
for $t>0$ and almost every $x$, where $a>0$. For every $\lambda>0$ and $u\in\mathcal D((A+\lambda)^{1/2})$,
$$
\|(A+\lambda)^{1/2}u\|_2^2\le1
\quad\Longrightarrow\quad
a\lambda\int_X(e^{|u|^2/a}-1)\dd\mu
\le B_{\sqrt2}+e-1.
$$
Neither finite total measure nor compact resolvent is required.
\end{proposition}

\begin{proof}
Set $B=A+\lambda$. Its kernel is $e^{-\lambda t}K_A(t,x,y)$. Cauchy--Schwarz against a kernel row gives
\begin{equation}\label{global:eq:smoothing}
\|e^{-tB/2}f\|_\infty
\le\sqrt{\frac at}\,e^{-\lambda t/2}\|f\|_2.
\end{equation}
The spectral theorem and Tonelli's theorem give
\begin{equation}\label{global:eq:frame}
\int_0^\infty\|Be^{-tB/2}u\|_2^2\dd t
=\|B^{1/2}u\|_2^2\le1,
\qquad
\lambda\|u\|_2^2\le1.
\end{equation}
Indeed, for a spectral value $\xi>0$,
$\int_0^\infty\xi^2e^{-t\xi}\dd t=\xi$.

For $\tau>0$, the scalar inequality $(1-e^{-x})^2\le x(1-e^{-x})$ implies
\begin{equation}\label{global:eq:high}
\|(I-e^{-\tau B})u\|_2^2
\le\tau\int_0^\tau\|Be^{-tB/2}u\|_2^2\dd t.
\end{equation}
The complementary spectral part is
$$
e^{-\tau B}u=\int_\tau^\infty Be^{-tB}u\dd t.
$$
This integral converges in $L^2$ because $B\ge\lambda>0$ and $\tau>0$. Applying \eqref{global:eq:smoothing} after factoring the semigroup at time $t/2$ gives
\begin{equation}\label{global:eq:low}
\|e^{-\tau B}u\|_\infty
\le\sqrt a\int_\tau^\infty
t^{-1/2}e^{-\lambda t/2}\|Be^{-tB/2}u\|_2\dd t.
\end{equation}
The right side is finite by \eqref{global:eq:frame} and
$\int_\tau^\infty t^{-1}e^{-\lambda t}\dd t<\infty$.
The measurable kernel gives an absolutely convergent pointwise integral almost everywhere. Testing on bounded functions of finite-measure support identifies it with the $L^2$ integral. Thus no assumption about strong measurability with values in $L^\infty$ is needed.

Let $u^*$ be the decreasing rearrangement of $|u|$ on $(0,\infty)$, extended by zero beyond $\mu(X)$ if the total measure is finite. Its positive level sets have finite measure since $u\in L^2$. A bounded term plus an $L^2$ term, followed by Chebyshev's inequality at the same measure $s$, gives
\begin{align}
u^*(s)\le{}&
\sqrt a\int_\tau^\infty
t^{-1/2}e^{-\lambda t/2}\|Be^{-tB/2}u\|_2\dd t
\notag\\
&+\sqrt{\frac{\tau}{s}}
\left(\int_0^\tau\|Be^{-tB/2}u\|_2^2\dd t\right)^{1/2}.
\label{global:eq:rearrange}
\end{align}
There is no change from $s$ to $s/2$: if the bounded summand has norm at most $c$, then the set $\{|u|>c+r\}$ is contained in the set where the remaining summand has modulus greater than $r$.

For $T\ge0$ choose
$$
s=\frac{e^{-T}}{a\lambda},\qquad
\tau=\frac{e^{-T}}{\lambda}=as,
$$
and define the nonnegative function on the whole real line
$$
\vartheta(S)=
\sqrt{\frac{e^{-S}}{\lambda}}\,
\left\|B e^{-e^{-S}B/(2\lambda)}u\right\|_2.
$$
The substitution $t=e^{-S}/\lambda$ in \eqref{global:eq:frame} yields
$\int_\R\vartheta^2\le1$. Divide \eqref{global:eq:rearrange} by $\sqrt a$. On $0<S<T$ discard the exponential weight, and on $S<0$ retain it. This gives
\begin{align}
\frac{u^*(e^{-T}/(a\lambda))}{\sqrt a}
\le{}&
\int_0^T\vartheta(S)\dd S
+\int_{-\infty}^0e^{-e^{-S}/2}\vartheta(S)\dd S
\notag\\
&+\left(\int_T^\infty\vartheta(S)^2\dd S\right)^{1/2}.
\label{global:eq:critical}
\end{align}
The coefficient of the first term is exactly $1$, which preserves the critical exponent. Moreover,
$$
\int_{-\infty}^0e^{-e^{-S}}\dd S
=\int_1^\infty\frac{e^{-r}}r\dd r\le e^{-1}<1.
$$
Cauchy--Schwarz, first on $(-\infty,0)$ and then in two components, bounds the last two terms in \eqref{global:eq:critical} by
$$
\sqrt2\left(
\int_{-\infty}^0\vartheta^2+
\int_T^\infty\vartheta^2\right)^{1/2}.
$$
Apply Lemma~\ref{two:lem:oned} to the restriction of $\vartheta$ to $[0,\infty)$, with $E_0=\int_{-\infty}^0\vartheta^2$ and $d=\sqrt2$. The change of variables $s=e^{-T}/(a\lambda)$ then gives
\begin{equation}\label{global:eq:largeramplitude}
\int_0^{1/(a\lambda)}e^{u^*(s)^2/a}\dd s
\le\frac{B_{\sqrt2}}{a\lambda}.
\end{equation}
This step controls the small-measure portion of the decreasing rearrangement, regardless of the total measure of $X$.

For $s\ge1/(a\lambda)$, \eqref{global:eq:frame} and Chebyshev imply
$$
0\le\frac{u^*(s)^2}{a}
\le\frac{\|u\|_2^2}{as}\le\frac1{a\lambda s}\le1.
$$
Using $e^z-1\le(e-1)z$ for $0\le z\le1$, we obtain
\begin{equation}\label{global:eq:smallamplitude}
\int_{1/(a\lambda)}^\infty
(e^{u^*(s)^2/a}-1)\dd s
\le\frac{e-1}{a}\|u\|_2^2
\le\frac{e-1}{a\lambda}.
\end{equation}
Equimeasurability, \eqref{global:eq:largeramplitude} and \eqref{global:eq:smallamplitude} prove the proposition.

All spectral identities hold for the quadratic form domain. Alternatively, begin with $1_{[\lambda,N]}(B)u$, which need not have finite rank. The energies do not increase, these truncations converge in the form norm, and an almost-everywhere subsequence followed by Fatou's lemma gives the same conclusion. No positivity preservation, Markov property or $L^1$ contraction is used in this abstract argument.
\end{proof}

\subsection{The two surface inequalities}

\begin{proof}[Proof of the inequalities in Theorem~\ref{thm:critical-exponential}]
On a relatively compact source domain $U$, Proposition~\ref{global:prop:heat} supplies the heat diagonal bound for $A_U=-4\Delta+|H|^2-b^2$ with $a=1/(16\pi)$. Its form domain is $W^{1,2}_0(U)$, or $W^{1,2}(\Sigma)$ on an entire compact source. Proposition~\ref{two:prop:abstract} proves part~(i).

For part~(ii), the same geometric proposition identifies the natural global form domain with
$$
\left\{v\in W^{1,2}(\Sigma):\int_\Sigma |H|^2v^2\dd\sigma<\infty\right\}.
$$
The energy bound is precisely $\|(A+\lambda)^{1/2}v\|_2^2\le1$. Proposition~\ref{global:prop:abstract}, again with $a=1/(16\pi)$, proves the global estimate. One common admissible constant for both parts is
$$
\mathcal C=\max\left\{B_{\sqrt{4e+1}},\,16\pi(B_{\sqrt2}+e-1)\right\}.
$$
It is independent of all geometric data and parameters.
\end{proof}

\subsection{Scaling and sharpness}

\subsubsection{Restoring arbitrary \texorpdfstring{$\lambda$}{lambda} by geometric scaling}

The proof above works for every $\lambda>0$ directly. It also agrees with the exact two-dimensional scaling. Set $g'=\lambda g$ on the ambient manifold and use the induced source metric. The connection and normal projections are unchanged, while
\begin{gather*}
d\sigma'=\lambda\,d\sigma,\qquad
|\nabla'v|_{g'_\Sigma}^2=\lambda^{-1}|\nabla v|_{g_\Sigma}^2,\\
H'=\lambda^{-1}H,\qquad
|H'|_{g'}^2=\lambda^{-1}|H|_g^2,\qquad
b'=\frac b{\sqrt\lambda}.
\end{gather*}
In particular $\Sec_{g'}\le-(b')^2$, and completeness and $C^{1,1}$ immersion are preserved. Both the Sobolev class and the finite $H$-weighted integral are preserved for fixed $\lambda$. Thus
\begin{align*}
&\int_\Sigma[4|\nabla'v|^2+(|H'|^2-(b')^2+1)v^2]\dd\sigma'\\
&\qquad=\int_\Sigma[4|\nabla v|^2+(|H|^2-b^2+\lambda)v^2]\dd\sigma,
\end{align*}
and
$$
\int_\Sigma(e^{16\pi v^2}-1)\dd\sigma'
=\lambda\int_\Sigma(e^{16\pi v^2}-1)\dd\sigma.
$$
Therefore the full-class $\lambda=1$ statement and Theorem~\ref{thm:critical-exponential}(ii) are equivalent. The curvature parameter is scaled together with the metric.

\subsubsection{The exponent}

Fix $k\ge1$ and $\lambda>0$. Take $\Sigma=\R^2$ as a linear plane in $\R^{2+k}$, with $b=0$ and $H=0$. For $T>1$ define
$$
w_T(r)=\frac1{2\sqrt{2\pi}}
\begin{cases}
\sqrt T,&0\le r\le e^{-T},\\
\log(1/r)/\sqrt T,&e^{-T}<r<1,\\
0,&r\ge1.
\end{cases}
$$
These continuous compactly supported functions belong to $W^{1,2}(\R^2)$, and direct polar integration gives
$$
4\int_{\R^2}|\nabla w_T|^2\dd x=1,\qquad
\|w_T\|_2^2
=\frac{T e^{-2T}}8+
\frac1{4T}\int_0^T s^2e^{-2s}\dd s
\sim\frac1{16T}.
$$
Set $v_T=w_T/\sqrt{1+\lambda\|w_T\|_2^2}$. The required energy is exactly $1$. For every $\alpha>16\pi$,
$$
\lambda\int_{\R^2}(e^{\alpha v_T^2}-1)\dd x
\ge\lambda\pi e^{-2T}
\left[
\exp\left(\frac{\alpha T}{8\pi(1+\lambda\|w_T\|_2^2)}\right)-1
\right]\longrightarrow\infty,
$$
since $\alpha/[8\pi(1+\lambda\|w_T\|_2^2)]\to\alpha/(8\pi)>2$. This proves optimality in part~(ii) even with $\lambda$ fixed. For part~(i), use $w_T$ without the $L^2$ normalization on $U=B_1$. Its gradient energy is one and $\int_U e^{\alpha w_T^2}\dd x\ge\pi\exp[(\alpha/(8\pi)-2)T]$, which diverges for every $\alpha>16\pi$.

\subsubsection{The global curvature coefficient}

Fix $\lambda>0$ and $\beta>1$. Choose one $b>0$ such that
$$
(\beta-1)b^2>\lambda.
$$
Let $\Sigma=\mathbb H^2(-b^2)$ be a totally geodesic plane in $\mathbb H^{2+k}(-b^2)$. This is a smooth complete embedded admissible source, and $H=0$.

For $\eta>1$ let
$$
w_\eta(r)=\cosh(br)^{-\eta/2}.
$$
This is smooth also at $r=0$, and polar integration, using the substitution $s=\cosh(br)$, gives
\begin{equation}\label{global:eq:hypmass}
\int_\Sigma w_\eta^2\dd\sigma=\frac{2\pi}{b^2(\eta-1)},
\qquad
\frac{4\int_\Sigma|\nabla w_\eta|^2\dd\sigma}
{\int_\Sigma w_\eta^2\dd\sigma}
=b^2\frac{2\eta^2}{\eta+1}.
\end{equation}
Thus $w_\eta\in W^{1,2}(\Sigma)$ and the displayed quotient tends to $b^2$ as $\eta\downarrow1$. Choose $\eta>1$ sufficiently close to $1$ that
$$
b^2\left(\frac{2\eta^2}{\eta+1}-\beta\right)+\lambda<0.
$$
Equation~\eqref{global:eq:hypmass} now gives
$$
\int_\Sigma[4|\nabla w_\eta|^2+(\lambda-\beta b^2)w_\eta^2]\dd\sigma<0.
$$
Fix this $b$, this geometry and this $\eta$. Every $v_T=Tw_\eta$, $T>0$, satisfies the modified energy bound, since its energy is negative. On any fixed small ball, $w_\eta$ has a positive lower bound. Consequently, for every positive exponent $\alpha$,
$$
\lambda\int_\Sigma(e^{\alpha T^2w_\eta^2}-1)\dd\sigma\longrightarrow\infty
\qquad(T\to\infty).
$$
Each individual exponential integral is finite, because $0<w_\eta\le1$ and
$e^{\alpha T^2w_\eta^2}-1\le(e^{\alpha T^2}-1)w_\eta^2$.
If compactly supported tests are desired, multiply the fixed $w_\eta$ by smooth radial cutoffs converging to one in $W^{1,2}$; its strictly negative modified energy persists for one sufficiently large fixed cutoff. Amplitude multiplication then gives the same contradiction.

This proves the asserted maximality across all $b\ge0$ for every prescribed $\lambda>0$. It does not assert failure for every fixed $b$; the inequality $(\beta-1)b^2>\lambda$ is essential to this particular counterexample.

\subsubsection{Finite-domain curvature sharpness}

For maximality of the curvature coefficient, fix $b=1$ and the totally geodesic plane $\mathbb H^2(-1)\subset\mathbb H^{2+k}(-1)$, for which $H=0$. For $a>1$, the smooth radial function $w_a(r)=\cosh(r)^{-a/2}$ satisfies
$$
\int_{\mathbb H^2}w_a^2\dd\sigma=\frac{2\pi}{a-1},
\qquad
\frac{4\int_{\mathbb H^2}|\nabla w_a|^2\dd\sigma}{\int_{\mathbb H^2}w_a^2\dd\sigma}
=\frac{2a^2}{a+1}\longrightarrow1
\quad(a\downarrow1).
$$
These identities follow by substituting $v=\cosh r$ into the radial integrals. Given $\beta>1$, choose $a>1$ with $2a^2/(a+1)<\beta$. Smooth radial cutoffs equal to one on $B_R$ and zero outside $B_{R+1}$, with uniformly bounded gradients, converge after multiplication by $w_a$ to $w_a$ in $W^{1,2}$. A sufficiently large fixed cutoff therefore gives a nonzero smooth compactly supported function $w$ with
$$
\int_{\mathbb H^2}(4|\nabla w|^2-\beta w^2)\dd\sigma<0.
$$
Fix a bounded smooth domain containing its support. Every multiple $sw$, $s>0$, has shifted energy at most one. For any $\alpha>0$, the integrals $\int_U e^{\alpha s^2w^2}\dd\sigma$ diverge as $s\to\infty$, since $|w|$ is bounded below by a positive number on some positive-area ball. This disproves the finite-domain assertion for every $\beta>1$, even at the fixed curvature $b=1$.

\section{The Nash inequalities and fixed-curvature sharpness}\label{sec:nash}

The inequalities below follow from Proposition~\ref{thm:sharp-sobolev} and Corollary~\ref{cor:homogeneous}. The equality argument in part~\textup{(i)} additionally determines the ambient geometry of the enclosed $(m+1)$-dimensional model ball.

\begin{corollary}\label{thm:nash}
Suppose $\Sec_M\le-b^2$ for some $b\ge0$.

\noindent\textup{(i)} If $m\ge3$, then for every $0\le\beta\le m^2$ and nonzero real $\varphi\in C_c^1(\Sigma)$,
\begin{align}\label{eq:nash}
\left(\int_\Sigma|\varphi|^2\dd\sigma\right)^{(m+2)/m}
\le{}& \frac{1}{m^2|\mathbb S^m|^{2/m}}
\int_\Sigma\left[\frac{4m}{m-2}|\nabla\varphi|^2
+(|H|^2-\beta b^2)|\varphi|^2\right]\dd\sigma\notag\\
&\hspace{18mm}\cdot\left(\int_\Sigma|\varphi|\dd\sigma\right)^{4/m}.
\end{align}
The constant is optimal over all the stated data for each fixed $k$ and is independent of $k$. Equality holds exactly when $(m^2-\beta)b^2=0$ and there is one connected component on which $\varphi$ is a nonzero constant and which is mapped, with one sheet, onto the boundary of an isometrically embedded totally geodesic $(m+1)$-dimensional ball of sectional curvature $-b^2$ in $M$. The function vanishes on all other components. The ball is Euclidean when $b=0$ and hyperbolic when $b>0$.

\noindent\textup{(ii)} If $m=2$, then for every real $\varphi\in C_c^1(\Sigma)$,
\begin{equation}\label{eq:nash-two}
\left(\int_\Sigma|\varphi|^2\dd\sigma\right)^2
\le\frac1{4\pi e}\int_\Sigma\bigl[4|\nabla\varphi|^2+(|H|^2-b^2)|\varphi|^2\bigr]\dd\sigma
\left(\int_\Sigma|\varphi|\dd\sigma\right)^2.
\end{equation}
No optimality or equality classification is asserted for the constant in \eqref{eq:nash-two}.
\end{corollary}

For $m\ge3$, the best constant over varying immersed submanifolds is attained by round spheres. This geometric optimization differs from the sharp Euclidean Nash inequality of Carlen and Loss \cite{CL93}, which concerns the Dirichlet energy on a fixed Euclidean space. At fixed $b>0$ and $\beta<m^2$, the same uniform constant is optimal but not attained; Corollary~\ref{cor:nash-range} gives the model sequence.

For the equality argument, we recall a consequence of stationary Euclidean monotonicity in the regularity class used here.

\begin{lemma}\label{lem:spherical}
A compact boundaryless $C^{1,1}$ weakly minimal immersion $J:K^m\to\mathbb S^{m+k}$ has source volume at least $|\mathbb S^m|$. If $K$ is connected, equality holds exactly for a one-sheet isometric parametrization of a great $m$-sphere.
\end{lemma}
\begin{proof}
The cone $(s,x)\mapsto sJ(x)$ is a stationary integral varifold in Euclidean space: minimality gives stationarity off the vertex, and the inner flux is bounded by $\varepsilon^m|K|\|X\|_\infty\to0$ for each compactly supported variation field $X$. Its normalized density at infinity is $|K|/|\mathbb S^m|$, while a nonvertex source sheet has density at least one. Stationary monotonicity \cite[Section~17]{Si83} gives the volume bound.

If $|K|=|\mathbb S^m|$, monotonicity is an equality at every nonvertex image point $z$. Its equality formula gives $(y-z)^\perp=0$ at almost every cone point $y$. Choose $m+1$ linearly independent image points $z_i$. Since the radial vector $y$ is tangent, all $z_i$ are tangent at almost every $y$, so the cone is contained in their fixed $(m+1)$-plane. Continuity extends the containment to the entire image. Thus $J$ is a local isometry onto the great sphere in that plane. Compactness gives a covering; simple connectivity of $\mathbb S^m$ and connectedness of $K$ give one sheet. Source multiplicities have been retained throughout.
\end{proof}

\begin{proof}[Proof of Corollary~\ref{thm:nash}(i)]
H\"older's inequality gives
$$
\left(\int|\varphi|^2\right)^{(m+2)/m}
\le\|\varphi\|_{2m/(m-2)}^2\left(\int|\varphi|\right)^{4/m}.
$$
For example, apply it to $|\varphi|^2=(|\varphi|^{2m/(m-2)})^{(m-2)/(m+2)}|\varphi|^{4/(m+2)}$ and raise to the power $(m+2)/m$. Proposition~\ref{thm:sharp-sobolev} and
$$
|H|^2-\beta b^2
=|H|^2-m^2b^2+(m^2-\beta)b^2
$$
prove \eqref{eq:nash} with the displayed constant. For $\Sigma=\mathbb S^m(R)\subset\mathbb R^{m+1}\subset\mathbb R^{m+k}$, with the latter inclusion linear, $b=0$ and a constant nonzero test, the Nash quotient is exactly $1/(m^2|\mathbb S^m|^{2/m})$. These are actual compactly supported tests on a compact source, so the uniform best value is attained.

Suppose equality holds. The added nonnegative term forces $(m^2-\beta)b^2=0$. Equality in H\"older forces $|\varphi|$ to take one positive value wherever it is nonzero. Continuity makes this positivity set a union of whole connected components, each compact because the support is compact. There is only one active component: if at least two critical-norm masses are positive, summing the componentwise Sobolev inequalities is strict because $(m-2)/m<1$. On that component $\varphi$ has a constant sign.

On the active component, $|\varphi|$ is a positive constant attaining equality in \eqref{eq:sharp-sobolev}. Use the probability $|\varphi|^{2m/(m-2)}d\sigma/\int|\varphi|^{2m/(m-2)}d\sigma$ to center the comparison map in Proposition~\ref{thm:sharp-sobolev}. The first and second variations along its coordinates, as in \cite[Section~3.3]{OAI26}, give
$$
\begin{aligned}
m^2|\mathbb S^m|^{2/m}\|\varphi\|_{2m/(m-2)}^2
&\le m\int |\varphi|^2\sum_{\alpha=1}^{m+k+1}|\nabla_\Sigma(F_\alpha\circ\iota)|^2\dd\sigma\\
&\le m^2\int\Omega^2|\varphi|^2\dd\sigma.
\end{aligned}
$$
Formula~\eqref{eq:vector-square} gives the reverse inequality. Hence $|\varphi|\Omega^{-(m-2)/2}$ is constant, the normal square vanishes, and equality holds in the tangent contraction. Thus $F\circ\iota$ is an isometric immersion for the source metric $\Omega^2g_\Sigma$. The normal equality says that the original immersion is weakly minimal for the conformal ambient metric, because
$$
H_{\Omega^2g}=\Omega^{-2}\bigl(H-m(\nabla^M\log\Omega)^\perp\bigr)=0.
$$
Minimality transfers to $F\circ\iota$ by first variation: the inequality $F^*g_{\mathbb S^{m+k}}\le\Omega^2g$ makes the difference of the two source volumes nonnegative under every compact variation, with value zero at the given immersion, so their first derivatives agree. The radial map is a diffeomorphism onto its image, and spherical variations near the compact image lift through it. The constant quotient is therefore
$$
m^2|\mathbb S^m|^{2/m}=m^2\left(\int\Omega^m\dd\sigma\right)^{2/m}.
$$
Lemma~\ref{lem:spherical} makes $F\circ\iota$ a one-sheet parametrization of a great $m$-sphere. For the radial parameters, the weights and last coordinates may be written as
$$
\begin{array}{lll}
b=0:&\displaystyle\Omega(r)=\frac{2t}{t^2+r^2},&\displaystyle F_{m+k+1}=\frac{r^2-t^2}{r^2+t^2},\quad t>0,\\[6pt]
b>0:&\displaystyle\Omega(r)=\frac{b\sqrt{1-\eta^2}}{\cosh(br)+\eta},&\displaystyle F_{m+k+1}=\frac{1+\eta\cosh(br)}{\cosh(br)+\eta},\quad -1<\eta<1.
\end{array}
$$
Since $\Omega$ is strictly decreasing in $r>0$, its constancy gives $r=R>0$. When $b=0$, the last centering coordinate gives $t=R$. When $b>0$, it gives $\eta=-1/\cosh(bR)$. In both cases $F$ has last coordinate zero on the source. Therefore the great sphere lies in that equator, and there is an $(m+1)$-dimensional linear subspace $W\subset T_pM$ such that the source is mapped once onto
$$
\{\exp_p(R\theta):\theta\in W,\ |\theta|=1\}.
$$

We verify the geometry of the enclosed $(m+1)$-dimensional ball. If $b=0$, the tangent contraction at $R$ gives $|Y_v(R)|=R|v|$ for every radial Jacobi field with initial data $Y_v(0)=0$, $Y_v'(0)=v$, where the initial radial direction $\theta$ belongs to $W$ and $v\in W\cap\theta^\perp$. Under $\Sec_M\le0$, the norm $|Y_v(s)|$ is convex, has initial derivative $|v|$, and is at least $s|v|$. Endpoint equality forces $|Y_v(s)|=s|v|$ on $[0,R]$. Polarization for all $v$ gives the Euclidean induced polar metric on
$$
D=\exp_p\{x\in W:|x|<R\}.
$$
This ball is totally geodesic in $M$. Indeed, the global radial comparison map $\log_p$ is $1$-Lipschitz, so for $x,y\in W$ with $|x|,|y|<R$,
$$
d_M(\exp_p x,\exp_p y)\ge|x-y|.
$$
The image of the straight segment from $x$ to $y$ stays in $D$ and has length $|x-y|$ by its induced metric. The opposite distance inequality follows. Thus straight segments map to ambient geodesics, proving total geodesicity and metric isometry of the ball, including its boundary.

If $b>0$, then $\beta=m^2$. At $r=R$ the parameter above gives $\Omega=b/\sinh(bR)$. Hence tangent contraction equality gives
$$
|Y_v(R)|=\frac{\sinh(bR)}{b}|v|,
$$
for all the same radial directions and transverse $v$ in $W$. Hessian comparison along the radial geodesic gives
$$
\frac{d}{ds}\log\frac{|Y_v(s)|}{\sinh(bs)/b}\ge0.
$$
The ratio tends to $|v|$ at zero and takes that value at $R$, so it is constant throughout. Polarization gives on $D$ the induced metric of the hyperbolic $(m+1)$-ball of curvature $-b^2$. To see total geodesicity, identify ambient polar coordinates with those of $\mathbb H^{m+k}(-b^2)$. Rauch makes this radial comparison map $1$-Lipschitz. For two points of $D$, their model geodesic stays inside the convex $(m+1)$-dimensional model ball determined by $W$. Its lift to $D$ has the same length. The comparison inequality and this lifted path give equality of ambient and model distances. Model geodesics therefore lift to ambient geodesics, proving that $D$ is isometrically embedded and totally geodesic.

The active source is the one-sheet boundary of $D$. Conversely, on the boundary of such a totally geodesic ball, $\varphi$ constant has the sharp Nash quotient whenever $(m^2-\beta)b^2=0$. The Euclidean calculation was given above. In the hyperbolic case the induced spherical radius is $\sinh(bR)/b$ and $|H|=mb\coth(bR)$, so
$$
|H|^2-m^2b^2=\frac{m^2}{(\sinh(bR)/b)^2},
$$
which gives the same sphere quotient. This proves the full equality classification.
\end{proof}

\begin{proof}[Proof of Corollary~\ref{thm:nash}(ii)]
It suffices to consider $\varphi\ne0$. By Lemma~\ref{lem:nonnegative}, the density $f=\varphi^2/\int_\Sigma\varphi^2\dd\sigma$ satisfies the homogeneous logarithmic inequality in dimension two. Jensen's inequality for the probability measure $f\dd\sigma$, applied on $\{f>0\}$, gives
$$
\int_\Sigma f\log f\dd\sigma\ge-2\log\int_\Sigma\sqrt f\dd\sigma
=\log\frac{\int_\Sigma\varphi^2\dd\sigma}{(\int_\Sigma|\varphi|\dd\sigma)^2}.
$$
Combining this with the case $m=2$, $\beta=1$ of \eqref{eq:homogeneous} and exponentiating proves \eqref{eq:nash-two}, since
$$
\int_\Sigma\left[\frac{|\nabla f|^2}{f}+(|H|^2-b^2)f\right]\dd\sigma
=\frac{\int_\Sigma[4|\nabla\varphi|^2+(|H|^2-b^2)\varphi^2]\dd\sigma}{\int_\Sigma\varphi^2\dd\sigma}.
$$
The compact support and local boundedness of $H$ ensure finite energy and absolute entropy. The weak gradient vanishes almost everywhere on the zero set, justifying the Fisher identity.
\end{proof}

\begin{corollary}\label{cor:nash-range}
Let $m\ge3$. For every fixed $b\ge0$, every fixed $0\le\beta\le m^2$ and every fixed $k$, the constant in Corollary~\ref{thm:nash} is optimal over the corresponding varying geometries. If $b>0$ and $\beta<m^2$, it is not attained, but shrinking hyperbolic geodesic spheres with constant tests form an extremizing sequence. For $\beta>m^2$ and every fixed $b>0$, no finite positive uniform Nash constant exists. Thus the maximal nonnegative curvature coefficient in a statement uniform over $b$ is $m^2$.
\end{corollary}

\begin{proof}
The upper bound is already proved. For $b=0$, Euclidean round spheres show sharpness for any $\beta$. Fix $b>0$ and take a round geodesic $m$-sphere of radius $R$ in a totally geodesic $\mathbb H^{m+1}(-b^2)\subset\mathbb H^{m+k}(-b^2)$. Its volume is $|\mathbb S^m|(\sinh(bR)/b)^m$ and its trace mean-curvature norm is $mb\coth(bR)$. The constant-test Nash quotient is exactly
$$
\frac{1}{|\mathbb S^m|^{2/m}\bigl[m^2+(m^2-\beta)b^2(\sinh(bR)/b)^2\bigr]}.
$$
For $\beta\le m^2$ this tends to $1/(m^2|\mathbb S^m|^{2/m})$ as $R\downarrow0$. For $\beta<m^2$ equality is excluded by Corollary~\ref{thm:nash}; for $\beta=m^2$ every $R$ already attains the best value. These spheres are genuine compact tests for each $R$, with no change of the fixed ambient curvature parameter. If $\beta>m^2$, then for sufficiently large $R$,
$$
|H|^2-\beta b^2=b^2\bigl(m^2\coth^2(bR)-\beta\bigr)<0.
$$
The same constant test has negative energy and positive Nash left side, proving failure in each fixed codimension.
\end{proof}

\textbf{Data Availability Statement}
Data sharing is not applicable to this article because no datasets were generated or analyzed during the current study.


\textbf{Conflict of interest} The authors declare that they have no conflict of interest.

\textbf{Code availability}
Not applicable.

\textbf{AI Assistance Statement}
 AI is used as an assistive tool.

\end{document}